\documentclass{amsart}
\usepackage{amsfonts}
\usepackage{pstricks,pstcol,pst-plot,pst-node,pst-tree,amssymb}
\usepackage{pstricks-add}
\usepackage{graphicx}

\numberwithin{equation}{section}
\input{tcilatex}

\begin{document}
\title[Flows and a tangency condition]{ Flows and a tangency condition for
embeddable $CR$ structures in dimension 3}
\author{ Jih-Hsin Cheng}
\address{Dedicated to Professor Stephen Yau \\
on his 60th birthday}
\address{Institute of Mathematics, Academia Sinica, Taipei, 11529 and
National Center for Theoretical Sciences, Taipei Office, Taiwan, R.O.C.}
\email{cheng@math.sinica.edu.tw}
\subjclass{2000 Mathematics Subject Classification. Primary 32G07, 32V30;
Secondary 32V20, 32V05. Research supported in part by National Science
Council grant NSC 101-2115-M-001-015-MY3 (R.O.C.).}

\begin{abstract}
We study the fillability (or embeddability) of 3-dimensional $CR$ structures
under the geometric flows. Suppose we can solve a certain second order
equation for the geometric quantity associated to the flow. Then we prove
that if the initial $CR$ structure is fillable, then it keeps having the
same property as long as the flow has a solution. We discuss the situation
for the torsion flow and the Cartan flow. In the second part, we show that
the above mentioned second order operator is used to express a tangency
condition for the space of all fillable or embeddable $CR$ structures at one
embedded in $\mathbb{C}^{2}.$
\end{abstract}

\maketitle


\section{Introduction}

A closed $CR$ manifold $M$ is fillable if $M$ bounds a complex manifold in
the smooth $(C^{\infty })$ sense (i.e. there exists a complex manifold with
smooth boundary $M$, and the complex structure restricts to the $CR$
structure on $M$). The notion of fillability is weaker than that of
embeddability. Recall that a $CR$ manifold is embeddable if it can be
embedded in $C^{N}$ for large $N$ with the $CR$ structure being the one
induced from the complex structure of $C^{N}$. The embeddability is a
special property for 3-dimensional $CR$ manifolds since any closed $CR$
manifold of dimension $\geq $ 5 is embeddable (\cite{BdM}). It is easy to
see that a closed embeddable (strongly pseudoconvex) $CR$ 3-manifold is
fillable by some well-known results (see the argument on page 543 in \cite%
{Ko}).

Conversely, if there exists a smooth strictly plurisubharmonic function
defined in a neighborhood of a fillable $M$, then $M$ is embeddable (\cite%
{Ko}, Theorem 5.3; in fact, any (compact) complex surface with nonempty
strongly pseudoconvex boundary can be made Stein by deforming it and blowing
down any exceptional curves according to \cite{Bo}).

In this paper, we first investigate the fillability of $3$-dimensional $CR$
structures under the (curvature) flows. Let $(M,\xi )$ denote a closed
(compact with no boundary) contact 3-manifold.with contact bundle
(structure) $\xi $ coorientable meaning that both $\xi $ and $TM/\xi ,$ the
normal bundle of $\xi $ in $TM,$ are orientable$.$ Let $J_{(t)}$ be a family
of $CR$ structures compatible with $\xi ,$ i.e., $J_{(t)}$ $\in $ $End(\xi )$
satisfying $J_{(t)}^{2}$ $=$ $-I.$ Set

\begin{equation}
\frac{\partial J_{(t)}}{\partial t}=2E_{J_{(t)}}  \label{1.1}
\end{equation}

\noindent where $E_{J_{(t)}}$ satisfies%
\begin{equation}
E_{J_{(t)}}\circ J_{(t)}+J_{(t)}\circ E_{J_{(t)}}=0.  \label{1.2}
\end{equation}

\noindent The most interesting case is to take $E_{J}$ to be the Cartan
(curvature) tensor $Q_{J}$ (see below for more details)$.$ Another
intriguing situation is to take $E_{J}$ to be the torsion tensor.

Choose a (global) contact form $\theta $ (exists since $TM/\xi $ is
orientable). We have a unitary coframe \{$\theta ^{1},$ $\theta ^{\bar{1}}\}$
and its dual $\{Z_{1},$ $Z_{\bar{1}}\}$ associated to $(J,\theta )$. Then we
have the torsion (tensor) $A_{11}$ and the Tanaka-Webster curvature $W,$
etc. (see Section 2 for details). We define a second order linear
differential operator $D_{J}$ from real functions to endomorphism fields by%
\begin{equation}
D_{J}f=(f_{,11}+iA_{11}f){\theta ^{1}\otimes }Z_{\bar{1}}+(f_{,\bar{1}\bar{1}%
}-iA_{\bar{1}\bar{1}}f){\theta ^{\bar{1}}\otimes }Z_{1}.  \label{1.3}
\end{equation}

\noindent We have lowered all the upper indices.

\bigskip

\textbf{Theorem A.} \textit{Suppose equation (\ref{1.1}) has a solution for }%
$0\leq t<\tau .$\textit{\ Suppose}%
\begin{equation}
J_{(t)}\circ D_{J_{(t)}}f+D_{J_{(t)}}g=E_{J_{(t)}}  \label{1.4}
\end{equation}%
\textit{\noindent has a solution of real functions (}$f,g$\textit{) defined
on }$M\times (0,\tau )$ \textit{with }$f$ $\neq $ $0.$\textit{\ Assume }$%
J_{(0)}$\textit{\ is fillable. Then }$J_{(t)}$\textit{\ is fillable for }$0$ 
$\leq $ $t$ $<$ $\tau .$

\bigskip

Recall that $D_{J_{(t)}}f=\frac{1}{2}L_{X_{f}}J_{(t)}$ (\cite{CL1}) in which 
$X_{f}=-fT$ $+$ $i(Z_{1(t)}f)Z_{\bar{1}(t)}$ $-i(Z_{\bar{1}(t)}f)Z_{1(t)}$
is the infinitesimal contact diffeomorphism induced by $f.$ So the image of $%
D_{J_{(t)}}$ describes the tangent space of the orbit of the symmetry group
acting on $J_{(t)}$ by the pullback (in this case, the contact
diffeomorphisms are our symmetries). Now condition (\ref{1.4}) means that $%
E_{J_{(t)}}$ sits in the \textquotedblright
complexification\textquotedblright\ of the infinitesimal orbit of contact
diffeomorphism group for all $t\in (0,\tau ).$

Write $u$ $=$ $f+ig$ and $E_{J_{(t)}}$ $=$ $E_{11(t)}\theta
_{(t)}^{1}\otimes Z_{\bar{1}(t)}$ $+$ $E_{\bar{1}\bar{1}(t)}\theta _{(t)}^{%
\bar{1}}\otimes Z_{1(t)}.$ We can then reformulate (\ref{1.4}) as an
equation for $u$:%
\begin{equation}
u_{,11}+iuA_{11(t)}=iE_{11(t)}  \label{1.5}
\end{equation}%
\textit{\noindent }(cf. (\ref{4.7}), (\ref{4.8}) in Section 4).\textit{\ }%
Note that equation (\ref{1.5}) has a solution $u$ $=$ $-1$ for $E_{11(t)}$ $%
= $ $-A_{11(t)}.$ Let $A_{J_{(t)}}$ or $A_{J_{(t)},\theta }$ denote the
torsion tensor $A_{11(t)}\theta _{(t)}^{1}\otimes Z_{\bar{1}(t)}$ $+$ $A_{%
\bar{1}\bar{1}(t)}\theta _{(t)}^{\bar{1}}\otimes Z_{1(t)}.$ We have the
following corollary.

\bigskip

\textbf{Corollary B.} \textit{Suppose the torsion flow}%
\begin{equation}
\frac{\partial J_{(t)}}{\partial t}=-2A_{J_{(t)}}  \label{1.6}
\end{equation}%
\textit{has a solution for }$0\leq t<\tau .$ \textit{Assume }$J_{(0)}$%
\textit{\ is fillable. Then }$J_{(t)}$\textit{\ is fillable for }$0\leq
t<\tau .$

\bigskip

Note that (\ref{1.6}) may not have a short-time solution for a general
smooth initial value. For the case $E_{J_{(t)}}$ $=$ $0,$ $J_{(t)}$ $=$ $%
J_{(0)}$ for all $t$ and $f$ $=$ $1,$ $g$ $=$ $0$ ($u$ $=$ $1,$ resp.$)$ is
a solution to (\ref{1.4}) ((\ref{1.5}), resp.) provided $A_{J_{(0)}}$ $=$ $0$
(and hence $A_{J_{(t)}}$ $=$ $0$ for all $t).$ In this case Theorem A is
obvious (note that $A_{J_{(0)}}$ $=$ $0$ implies that $J_{(0)}$ is
embeddable, and hence fillable. See the remark after Corollary C). For the
case $E_{J_{(t)}}$ $=$ $Q_{J_{(t)}},$ we do have a short-time solution for (%
\ref{1.1}) (see \cite{CL1}). In \cite{CL1}, we study an evolution equation
for $CR$ structures $J_{(t)}$ on $(M,\xi )$ according to their Cartan
(curvature) tensor $Q_{J_{(t)}}$(see also Section 2):%
\begin{equation}
\frac{\partial J_{(t)}}{\partial t}=2Q_{J_{(t)}}.  \label{1.7}
\end{equation}

\noindent We will often call this evolution equation (\ref{1.7}) the Cartan
flow. Since (\ref{1.7}) is invariant under a big symmetry group, namely, the
contact diffeomorphisms, we add a gauge-fixing term on the right-hand side
to break the symmetry. The gauge-fixed (called \textquotedblright
regularized\textquotedblright\ in \cite{CL1}) Cartan flow reads as follows:%
\begin{equation}
\frac{\partial J_{(t)}}{\partial t}=2Q_{J_{(t)}}-\frac{1}{6}%
D_{J_{(t)}}F_{J_{(t)}}K  \label{1.8}
\end{equation}

\noindent (see \cite{CL1} or Section 2 for the meaning of notations). Now it
is natural to ask the following question:

\medskip

\textbf{Question 1.1}: \textit{Is the fillability (or embeddability)
preserved under the (gauge-fixed) Cartan flow (\ref{1.7}) (or (\ref{1.8}))?}

\bigskip

An affirmative answer to the above question has an application in
determining the topology of the space of all fillable $CR$ structures. For
instance, one can apply such a result plus the convergence of the long time
solution to (\ref{1.8}) (expected for $S^{3}$) to prove that the space of
all fillable $CR$ structures on $S^{3}$ is contractible (cf. Remark 4.3 in 
\cite{El}). For other topological applications of solving (\ref{1.8}), we
refer the reader to \cite{Ch}.

In view of Theorem A we make the following conjecture.

\bigskip

\textbf{Conjecture 1.2}: \textit{We can find real functions }$f\neq 0$%
\textit{\ and }$g$\textit{\ (}$u$\textit{\ with }$\func{Re}u$\textit{\ }$%
\neq $\textit{\ }$0,$\textit{\ resp.) such that }$J\circ D_{J}f+D_{J}g=Q_{J}$%
\textit{\ (}$u_{,11}+iuA_{11}=iQ_{11},$\textit{\ resp.) for }$J$\textit{\
fillable or embeddable.}

\bigskip

On the other hand, we examine a family of $CR$ structures embedded in $%
\mathbb{C}^{2}$ and find its tangent to be of the form $J\circ D_{J}f+D_{J}g$
(see Section 6)$.$ Let $W_{J,\theta }$ denote the Tanaka-Webster curvature
of a pseudohermitian structure $(J,\theta )$ (see Section 2 for the
definition).

\bigskip

\textbf{Corollary C.} \textit{Suppose }$J_{(0)}$\textit{\ is fillable with }$%
A_{J_{(0),}\theta }=0$\textit{\ and }$W_{J_{(0)},\theta }>0$\textit{\ (or }$%
<0$\textit{, resp.). Then the solution }$J_{(t)}$\textit{\ to (\ref{1.8})
with }$K=J_{(0)}$\textit{\ stays fillable for a short time.}

\bigskip

The idea of the proof of Corollary C is to show $A_{J_{(t),}\theta }=0$ for
a short time (see Lemma 3.1) and then solve (\ref{1.4}) for $E_{J_{(t)}}$ $=$
$Q_{J_{(t)}}-\frac{1}{12}D_{J_{(t)}}F_{J_{(t)}}K$ in Theorem A (see Section
3).

The proof of Theorem A\textit{\ }is a direct construction of an integrable
almost complex structure $\widetilde{J}$ on $M\times (0,\tau )$ so that $%
\widetilde{J}|_{\xi }=J_{(t)}$ at $M\times \{t\}$ (see Section 4 for
details). Then we glue this complex structure $\widetilde{J}$ with the one
induced by the complex surface that $(M,J_{(0)})$ bounds along $M\times
\{0\} $ (identified with $M$). After we obtained the above result, L\'{a}szl%
\'{o} Lempert pointed out to the author that the existence of a $CR$ vector
field $T$ is sufficient to imply the embeddability of the $CR$ structure
(see \cite{Lem}). So by Lemma 3.1\textit{\ }we can remove the condition\ in
Corollary\textit{\ C} on the Tanaka-Webster curvature according to \cite{Lem}%
. We speculate that the embeddability (or fillability) is preserved under
the (gauge-fixed) Cartan flow without any conditions (Question 1.1).

Our method can also be applied to the problem of local embeddability. An
obvious case is that $(M,J_{(0)})$ with $A_{J_{(0)},\theta }$ $=$ $0$ is
locally embeddable. The reason is that the above $\widetilde{J}$ is
integrable on $M\times (-\tau ,\tau )$ since the torsion flow (\ref{1.6})
has an obvious solution $J_{(t)}$ $=$ $J_{(0)}$ for $t$ $\in $ $(-\tau ,\tau
)$ (with $A_{J_{(t)}}$ $=$ $A_{J_{(t)},\theta }$ $=$ $0).$ We single it out
as a corollary.

\bigskip

\textbf{Corollary D.} $(M,J_{(0)})$\textit{\ with }$A_{J_{(0)},\theta }$%
\textit{\ }$=$\textit{\ }$0$\textit{\ is locally embeddable.}

\bigskip

We will further study the local embedding problem in a forthcoming paper.

Another remark is that we may couple equation (\ref{1.1}) with an evolution
equation in contact form $\theta :$%
\begin{equation}
\frac{\partial \theta _{(t)}}{\partial t}=2h\theta _{(t)}.  \label{1.8.1}
\end{equation}

\noindent This does not affect the condition (\ref{1.4}) for integrability
of $\widetilde{J}.$ For instance, we may take $h$ $=$ $W_{J_{(t)},\theta
_{(t)}}$ to couple with the (pure) torsion flow (\ref{1.6}). (this coupled
torsion flow is the negative gradient flow of the pseudohermitian
Einstein-Hilbert action : $-\int W_{J,\theta }\theta \wedge d\theta $). See
a recent paper \cite{CKW} for more information about this flow.

One may speculate that $J\circ D_{J}f+D_{J}g$ is a tangent at $J$ of the
space of all embeddable or fillable (compatible) $CR$ structures on a fixed
contact 3-manifold. Starting from Section 5, we will justify this statement
for $J$ associated to a real hypersurface embedded in $\mathbb{C}^{2}$.

As we mentioned in the beginning, the (global) embeddability of a compact
(strongly pseudoconvex) $CR$ manifold (of hypersurface type) is of special
interest in dimension 3 since it is always embeddable for dimension $\geq 5$
and not always so for dimension $=3$. The analytic reason is that the
operator $\bar{\partial}_{b}$ associated to the concerned $CR$ structure is
solvable for type (0,1)-form when dimension $\geq 5$ and the space is
compact (\cite{FK}, \cite{BdM}). In dimension 3, the analysis of $\bar{%
\partial}_{b}$ is delicate. Lewy's example (\cite{Lew}) tells that there are
no solutions at all for the equation $\bar{\partial}_{b}u=\psi $ even with
certain $C^{\infty }$ functions $\psi $ (or, say, (0,1)-form). Using $\psi $%
, Nirenberg (\cite{N}) was able to find an example which is not embeddable
locally. On the other hand, an example of real-analytic perturbation of $%
S^{3}$ was constructed (\cite{R}, \cite{AS}), which (surely is locally
embeddable) is not globally embeddable. It is now understood that $CR$
structures with \textquotedblleft exotic" underlying contact structures are
generally non-embeddable. (\cite{El-1}) In \cite{BE}, Burns and Epstein made
a detailed study on perturbations of the standard $CR$ structure on the
three sphere. They gave a \textquotedblleft pointwise" criterion for
embeddability in terms of the spectrum of $\square _{b}$. Also they showed
that structures which are infinitesimally obstructed can not be embedded as
\textquotedblleft small perturbations" of the standard sphere. In fact, an
explicit decomposition of the tangent space $T_{0}$ to the perturbations of $%
S^{3}$ is given in their paper: 
\begin{equation*}
T_{0}=N\oplus E\oplus O
\end{equation*}%
where $O$ is tangent to the orbit of symmetries (in this case, they are
contact diffeomorphisms of $S^{3}$), $E$ is tangent to a family of
perturbations all of which are embedded in $\mathbb{C}^{2}$ as small
perturbations of $S^{3}$, and $N$ is tangent to a family of structures which
generically embed in no $\mathbb{C}^{n}$ , and none of which can embed in $%
\mathbb{C}^{2}$ as small perturbations of $S^{3}$. Later J. Bland, T.
Duchamp and L. Lempert studied the case of $CR$ structures induced by
strictly linearly convex domains extending the above case.

In general, we want to give a qualitative description of a11 embeddable $CR$
structures near an embedded one in smooth tame or Banach (if possible)
category. Let $(M,\xi )$ denote a closed (compact with no boundary) contact
3-manifold.with contact structure $\xi $ coorientable. The set of all $CR$
structures compatible with $\xi $ is denoted by $\mathfrak{J}_{\xi }$. It is
known \cite{CL1} that $\mathfrak{J}_{\xi }$ can be parametrized by sections
of a certain real 2-dimensional subbundle of $End(\xi )$. Therefore $%
\mathfrak{J}_{\xi }$ is a tame Fr\'{e}chet or Banach manifold. Our goal is
to solve the following problem:

\bigskip

\textbf{Conjecture 1.3.} \textit{Given a CR structure }$\mathit{(M,\xi ,J)}$ 
\textit{embedded in $\mathbb{C}^{2}$, there is a submanifold $\mathfrak{J}%
_{c}\subset \mathfrak{J}_{\xi }$ passing through $J$ and a neighborhood $%
\mathfrak{U}$ of $J$ in $\mathfrak{J}_{\xi }$ such that $\tilde{J}\in 
\mathfrak{U}$ is embeddable in $\mathbb{C}^{2}$ and realized by a nearby
embedding if and only if $\tilde{J}\in \mathfrak{J}_{c}$.}

\bigskip

In this paper, we look at the above problem infinitesimally. We will
describe the tangent space $T_{J}\mathfrak{J}_{c}$, of $\mathfrak{J}_{c}$,
at $J$ (meaning the space of all $\frac{d}{dt}|_{t=0}$ $J_{(t)}$ for $%
J_{(t)} $ $\in $ $\mathfrak{J}_{c}$ with $J_{(0)}$ $=$ $J)$ as $J\circ
D_{J}f+D_{J}g$ for all real functions $f$ and $g$. The statment
\textquotedblleft and realized by a nearby embedding" in Conjecture 1.3 may
be removed due to a result of Lempert \cite{Lem} for the stability of
embeddings if $(M,\xi ,J)$ is the boundary of a strictly linearly convex
domain in $\mathbb{C}^{2}$. But in general unstable $CR$ embeddings do exist
(\cite{CaL}).

Let $M\subset \mathbb{C}^{2}$ be a closed strongly pseudoconvex real
hypersurface with induced $CR$ structure $(\xi ,J)$. Let $\mathfrak{E}_{c}$
denote the set of all contact embeddings: ($M,\xi )$ $\rightarrow $ ${%
\mathbb{C}^{2}}$ near the inclusion map (see Section 5). Two contact
embeddings $\varphi ,$ $\psi $ are equivalent (in notation, $\varphi \sim
\psi )$ if $\varphi ^{\ast }J_{\mathbb{C}^{2}}=\psi ^{\ast }J_{\mathbb{C}%
^{2}}.$ Define (w.r.t. unitary (co)frame)

\begin{equation*}
\mathfrak{D}_{J}h=(h_{,\bar{1}\bar{1}}-iA_{\bar{1}\bar{1}}h)\theta ^{\bar{1}%
}\otimes Z_{1}
\end{equation*}

\noindent for $h$ $\in $ $C^{\infty }(M,\mathbb{C}),$ say (cf. (\ref{6.9}).
Note that 
\begin{equation*}
2\func{Re}\mathfrak{D}_{J}h=J\circ D_{J}f+D_{J}g
\end{equation*}%
\noindent for $h$ $=$ $g+if.$ Define a type (1,0) vector field $Y_{h}$ by%
\begin{equation*}
Y_{h}:=i(h\zeta +h^{,1}Z_{1})
\end{equation*}%
\noindent where $i\zeta ,$ $Z_{1}$ are type (1,0) (local) vector fields dual
to $\theta $ $=$ $-i\partial \gamma ,$ $\theta ^{1}$ near $M$ $:=$ $\{\gamma 
$ $=$ $0\}$ in $C^{2}$ (see Section 6)$.$

\bigskip

\textbf{Theorem D. }\textit{In the situation described above, we have the
following commutative diagram:}%
\begin{equation}
\begin{array}{ccccc}
\hbox{(complex version)}T_{[i_{M}]}(\mathfrak{E}_{c}/_{\sim }) & \overset{%
4Im\circ \bar\partial_b}{\longrightarrow} & T_{J}\mathfrak{J}_{c} & \subset
& T_{J}\mathfrak{J}_{\xi } \\[8pt] 
\lbrack Y_{\cdot }]\uparrow &  & 4\mathbb{R}e\uparrow &  & 4\mathbb{R}%
e\uparrow \\[8pt] 
C^{\infty }(M,\mathbb{C})/Ker\mathfrak{D}_{J} & \overset{\mathfrak{D}_{J}}{%
\longrightarrow} & \mathrm{Range}\ \mathfrak{D}_{J} & \subset & \Gamma
(T_{0,1}^{\ast }(M)\otimes T_{1,0}(M))%
\end{array}
\label{1.9}
\end{equation}

\noindent \textit{where the maps indicated by "}$\rightarrow "$\textit{\ or "%
}$\uparrow "$ \textit{are all one-one correspondences.}

\bigskip

The results in this paper were obtained in nineties. The author would like
to thank L\'{a}szl\'{o} Lempert, Jack Lee, and I-Hsun Tsai for discussions
during those years. In recent years, embeddability has played an important
role in the study of $CR$ positive mass theorem of $3D$ and positivity of $CR
$ Paneitz operator (\cite{CMY}, \cite{CCY1}, \cite{CCY2}). In this respect,
the author would like to thank Paul Yang, Andrea Malchiodi, and Hung-Lin
Chiu for many useful discussions. Theorem D is related to such a study. 

\bigskip\ 

\section{Review in $CR$ and pseudohermitian geometries}

For most of basic material we refer the reader to \cite{We}, \cite{Ta} or 
\cite{L1}. Throughout the paper, our base space $M$ is a closed (compact
with no boundary) contact 3-manifold with a cooriented contact structure $%
\xi $ meaning that $\xi $ and $TM/\xi $ are orientable. A $CR$ structure $J$
(compatible with $\xi $) is an endomorphism on $\xi $ with $J^{2}=-identity$.

By choosing a (global) contact form $\theta $ (exists since the normal
bundle $TM/\xi $ of $\xi $ in $TM$ is orientable), we can talk about
pseudohermitian geometry. The Reeb vector field $T$ is uniquely determined
by ${\theta }(T)=1$ and $T{\rfloor }d{\theta }=0.$ We choose a (local)
complex vector field $Z_{1}$, an eigenvector of $J$ with eigenvalue $i$, and
a (local) complex 1-form ${\theta }^{1}$ such that $\{{\theta },{\theta ^{1}}%
,{\theta ^{\bar{1}}}\}$ is dual to $\{T,Z_{1},Z_{\bar{1}}\}$ (here ${\theta
^{\bar{1}}}$ and $Z_{\bar{1}}$ mean the complex conjugates of ${\theta }^{1}$
and $Z_{1}$ resp.). It follows that $d{\theta }=ih_{1{\bar{1}}}{\theta ^{1}}{%
\wedge }{\theta ^{\bar{1}}}$ for some nonzero real function $h_{1{\bar{1}}}$
(may assume $h_{1{\bar{1}}}$ $>$ $0;$ otherwise, replace $\theta $ by $%
-\theta )$. We can then choose a $Z_{1}$ (hence $\theta ^{1}$) such that $%
h_{1{\bar{1}}}=1$. That is to say%
\begin{equation}
d{\theta }=i{\theta ^{1}}{\wedge }{\theta ^{\bar{1}}}.  \label{2.1}
\end{equation}

We will always assume our pseudohermitian structure $(J,{\theta })$
satisfies (\ref{2.1}), i.e., $h_{1{\bar{1}}}=1$ throughout the paper. The
pseudohermitian connection of $(J,{\theta })$ is the connection ${\nabla }^{{%
\psi }.h.}$ on $TM{\otimes }C$ (and extended to tensors) given by

\begin{equation*}
{\nabla }^{{\psi }.h.}Z_{1}={\omega _{1}}^{1}{\otimes }Z_{1},{\nabla }^{{%
\psi }.h.}Z_{\bar{1}}={\omega _{\bar{1}}}^{\bar{1}}{\otimes }Z_{\bar{1}},{%
\nabla }^{{\psi }.h.}T=0
\end{equation*}

\noindent in which the connection 1-form ${\omega _{1}}^{1}$ is uniquely
determined by the following equation and associated normalization condition:%
\begin{eqnarray}
d{\theta ^{1}} &=&{\theta ^{1}}{\wedge }{\omega _{1}}^{1}+{A^{1}}_{\bar{1}}{%
\theta }{\wedge }{\theta ^{\bar{1}},}  \label{2.2} \\
{0} &=&{\omega _{1}}^{1}+{\omega _{\bar{1}}}^{\bar{1}}.  \notag
\end{eqnarray}

\noindent The coefficient ${A^{1}}_{\bar{1}}$ in (\ref{2.2}) and its complex
conjugate ${A^{\bar{1}}}_{1}$ are components of the torsion (tensor) $%
A_{J,\theta }$ $=$ ${A^{\bar{1}}}_{1}{\theta ^{1}\otimes }Z_{\bar{1}}+{A^{1}}%
_{\bar{1}}{\theta ^{\bar{1}}\otimes }Z_{1}.$ Since $h_{1{\bar{1}}}=1$, $A_{{%
\bar{1}}{\bar{1}}}=h_{1{\bar{1}}}{A^{1}}_{\bar{1}}={A^{1}}_{\bar{1}}$.
Further $A_{11}$ is just the complex conjugate of $A_{{\bar{1}}{\bar{1}}}$.\
Write $J=i$ ${\theta ^{1}\otimes }Z_{1}-i{\theta ^{\bar{1}}\otimes }Z_{\bar{1%
}}$. It is not hard to see from (\ref{2.1}) and (\ref{2.2}) that%
\begin{equation}
L_{T}J=2J\circ A_{J,\theta }  \label{2.3}
\end{equation}

\noindent where $L_{T}$ denotes the Lie differentiation in the direction $T$
(this is the case when $f=-1$ in Lemma 3.5 of \cite{CL1}). So the vanishing
torsion is equivalent to $T$ being an infinitesimal $CR$ diffeomorphism. We
can define the covariant differentiations with respect to the
pseudohermitian connection. For instance, $f_{,1}=Z_{1}f$, $f_{,1{\bar{1}}%
}=Z_{\bar{1}}Z_{1}f-{\omega _{1}}^{1}(Z_{\bar{1}})Z_{1}f$ for a (smooth)
function $f$ (see, e.g., Section 4 in \cite{L1}). Now differentiating ${%
\omega _{1}}^{1}$ gives%
\begin{equation}
d{\omega _{1}}^{1}=W{\theta ^{1}}{\wedge }{\theta ^{\bar{1}}}+2iIm(A_{11,{%
\bar{1}}}{\theta ^{1}}{\wedge }{\theta })  \label{2.4}
\end{equation}

\noindent where $W$ or $W_{J,\theta }$ (to emphasize the dependence of the
pseudohermitian structure) is called the (scalar) Tanaka-Webster curvature.

There are distinguished $CR$ structures $J$, called spherical, if $(M,\xi
,J) $ is locally $CR$ equivalent to the standard 3-sphere $(S^{3},\hat{\xi},%
\hat{J})$, or equivalently if there are contact coordinate maps into open
sets of $(S^{3},\hat{\xi})$ so that the transition contact maps can be
extended to holomorphic transformations of open sets in $\mathbb{C}^{2}$. In
1930's, Elie Cartan (\cite{Ca}; see also \cite{CL1}) obtained a geometric
quantity, denoted as $Q_{J}$, by solving the local equivalence problem for $%
3 $-dimensional $CR$ structures so that the vanishing of $Q_{J}$
characterizes $J$ to be spherical. We will call $Q_{J}$ the Cartan
(curvature) tensor. Note that $Q_{J}$ depends on a choice of contact form $%
\theta .$ It is $CR$-covariant in the sense that if $\tilde{\theta}%
=e^{2f}\theta $ is another contact form and $\tilde{Q}_{J}$ is the
corresponding Cartan tensor, then $\tilde{Q}_{J}=e^{-4f}Q_{J}$. We can
express $Q_{J}$ in terms of pseudohermitian invariants. Write $Q_{J}=iQ_{11}{%
\theta ^{1}\otimes }Z_{\bar{1}}-iQ_{\bar{1}\bar{1}}{\theta ^{\bar{1}}\otimes 
}Z_{1}$ (note that $Q_{1}{}^{\bar{1}}=Q_{11}$ and $Q_{\bar{1}}{}^{1}=Q_{\bar{%
1}\bar{1}}$ since we always assume $h_{1{\bar{1}}}=1$). We have the
following formula (Lemma 2.2 in \cite{CL1}):%
\begin{equation}
Q_{11}=\frac{1}{6}W_{,11}+\frac{i}{2}WA_{11}-A_{11,0}-\frac{2i}{3}A_{11,\bar{%
1}1.}  \label{2.5}
\end{equation}

In terms of local coframe fields we can express the Cartan flow (\ref{1.7})
as follows:%
\begin{equation}
\ \dot{\theta}^{1}=-Q_{\bar{1}\bar{1}}\theta ^{\bar{1}}  \label{2.6}
\end{equation}

\noindent (cf. (2.16) in \cite{CL1} with $E_{\bar{1}}$ $^{1}$ replaced by $%
-iQ_{\bar{1}\bar{1}}$). The torsion evolves under the Cartan flow as shown
in the follow formula:%
\begin{equation}
\dot{A}_{11}=-Q_{11,0}  \label{2.7}
\end{equation}

\noindent (this is the complex conjugate of (2.18) in \cite{CL1} with $E_{%
\bar{1}}$ $^{1}$ replaced by $-iQ_{\bar{1}\bar{1}}$). Since the Cartan flow
is invariant under the pullback action of contact diffeomorphisms (cf. the
argument in the proof of Proposition 3.6 in \cite{CL1}), we need to add a
gauge-fixing term to the right-hand side of (\ref{1.7}) to get the
subellipticity of its linearized operator. Let us recall what this term is.
First we define a quadratic differential operator $F_{J}$ from endomorphism
fields to functions by ( \cite{CL1}, p.236 and note that $h_{1\bar{1}}=1$
here)%
\begin{equation}
F_{J}E=(iE_{1\bar{1}}E_{11,\bar{1}\bar{1}}+iE_{\bar{1}\bar{1}}E_{11,\bar{1}%
1})+conjugate.  \label{2.8}
\end{equation}

Also we define a linear differential operator $D_{J}$ from functions to
endomorphism fields and its formal adjoint $D_{J}^{\ast }$ by%
\begin{eqnarray}
D_{J}f &=&(f_{,11}+iA_{11}f){\theta ^{1}\otimes }Z_{\bar{1}}+(f_{,\bar{1}%
\bar{1}}-iA_{\bar{1}\bar{1}}f){\theta ^{\bar{1}}\otimes }Z_{1},  \label{2.9}
\\
\ D_{J}^{\ast }E &=&E_{11,\bar{1}\bar{1}}+E_{\bar{1}\bar{1},11}-iA_{\bar{1}%
\bar{1}}E_{11}+iA_{11}E_{\bar{1}\bar{1}}  \notag
\end{eqnarray}

\noindent (note that we have used the notations $D_{J},D_{J}^{\ast }$
instead of $B_{J}^{^{\prime }},B_{J}$ in \cite{CL1}, resp.). Now let $K$ be
a fixed $CR$ structure. The Cartan flow with a gauge-fixing term reads as
follows: (this is (\ref{1.8}))%
\begin{equation*}
\frac{\partial J_{(t)}}{\partial t}=2Q_{J_{(t)}}-\frac{1}{6}%
D_{J_{(t)}}F_{J_{(t)}}K.
\end{equation*}

We also need the following commutation relations often:%
\begin{eqnarray}
C_{I,01}-C_{I,10} &=&C_{I,\bar{1}}A_{11}-kC_{I}A_{11,\bar{1}}  \label{2.10}
\\
C_{I,0\bar{1}}-C_{I,\bar{1}0} &=&C_{I,1}A_{\bar{1}\bar{1}}+kC_{I}A_{\bar{1}%
\bar{1},1}  \notag \\
C_{I,1\bar{1}}-C_{I,\bar{1}1} &=&iC_{I,0}+kC_{I}W.  \notag
\end{eqnarray}

\noindent Here $C_{I}$ denotes a coefficient of a tensor with multi-index $I$
consisting of $1$ and $\bar{1},$ and $k$ is the number of $1$ in $I$ minus
the number of $\bar{1}$ in $I$ (an extension of formulas in \cite{L2}).

\bigskip

\section{Proof of Corollary C}

\textbf{Lemma 3.1}\textit{. Suppose there is a contact form }$\theta $%
\textit{\ such that the torsion }$A_{J_{(0),}\theta }$\textit{\ vanishes.
Then under the gauge-fixed Cartan flow (\ref{1.8}) (assuming smooth
solution) with }$K=J_{(0)}$\textit{, }$A_{J_{(t),}\theta }$\textit{\ stays
vanishing.}

\bigskip

Let $T$ denote the Reeb vector field associated with the contact form $%
\theta $. The vanishing of the torsion is equivalent to saying that $T$ \ is
an infinitesimal $CR$ diffeomorphism (see (\ref{2.3})). We say a $CR$
manifold has transverse symmetry if the infinitesimal generator of a
one-parameter group of $CR$ diffeomorphisms is everywhere transverse to $\xi 
$. Such an infinitesimal generator can be realized as the Reeb vector field
for a certain contact form $\hat{\theta}$ (\cite{L2}). From Lemma 3.1 we have

\bigskip

\textbf{Corollary 3.2.}\textit{\ The }$CR$\textit{\ structures }$J_{(t)}$%
\textit{\ stay having the same transverse symmetry as }$J_{(0)}$\textit{\
does under the gauge-fixed Cartan flow (\ref{1.8})} \textit{with }$K=J_{(0)}$%
\textit{\ and }$\theta =\hat{\theta}.$

\bigskip

\proof
\textbf{(of Lemma 3.1)} We will compute the evolution of the torsion under
the flow (\ref{1.8}) (with $K$ being the initial $CR$ structure $J_{(0)}$).
First, instead of (\ref{2.7}), we have%
\begin{equation}
\dot{A}_{11}=-Q_{11,0}-\frac{i}{12}(D_{J}F_{J}K)_{11,0}.  \label{3.1}
\end{equation}

From the formula (\ref{2.5}) for $Q_{11}$, we compute $Q_{11,0}$. Using the
commutation relations (\ref{2.10}) and the Bianchi identity: $W_{,0}=A_{11,%
\bar{1}\bar{1}}+A_{\bar{1}\bar{1},11}$ (\cite{L2}), we can express $Q_{11,0}$
only in terms of $A_{11,}A_{\bar{1}\bar{1}}$ and their covariant derivatives
as follows:%
\begin{equation}
Q_{11,0}=\frac{1}{6}(A_{11,\bar{1}\bar{1}11}+A_{\bar{1}\bar{1}%
,1111})-A_{11,00}-\frac{2i}{3}A_{11,\bar{1}10}+l.w.t..  \label{3.2}
\end{equation}

\noindent where $l.w.t.$ means a lower weight term in $A_{11}$ and $A_{\bar{1%
}\bar{1}}.$We count covariant derivatives in $1$ or $\bar{1}$ direction ($0$
direction, resp.) as weight $1$ (weight $2$, resp.) and we call a term of
weight $m$ if its total weight of covariant derivatives is $m.$ For
instance, $A_{11,\bar{1}\bar{1}11}$, $A_{11,00}$ and $A_{11,\bar{1}10}$ are
all of weight 4. So more precisely each single term in $l.w.t.$ must contain
terms of weight $\leq 3$ in $A_{11}$ or $A_{\bar{1}\bar{1}}.$In particular,
if $A_{11}=0$, then $l.w.t.=0.$ Note that $A_{\bar{1}\bar{1},1111}$ is a
\textquotedblright bad\textquotedblright term in the sense that we need a
gauge-fixing term to cancel it and obtain a fourth order subelliptic
operator in $A_{11}.$ Now by (\ref{2.9}) the gauge-fixing term in (\ref{3.1}%
) (up to a multiple) reads as%
\begin{eqnarray}
(D_{J}F_{J}K)_{11,0} &=&(F_{J}K)_{,110}+i[A_{11}(F_{J}K)]_{,0}  \label{3.3}
\\
&=&(F_{J}K)_{,011}+l.w.t.(\text{in }A_{11})  \notag
\end{eqnarray}

\noindent (we have used the commutation relations (\ref{2.10}) for the last
equality). Write $K=K_{11}\theta ^{1}\otimes Z_{\bar{1}}+$ $K_{1\bar{1}}$\ ${%
\theta ^{1}\otimes }Z_{1}+K_{\bar{1}\bar{1}}$\ ${\theta ^{\bar{1}}\otimes }%
Z_{1}+K_{\bar{1}1}{\theta ^{\bar{1}}\otimes }Z_{\bar{1}}$ where $K_{\bar{1}%
\bar{1}},K_{\bar{1}1}$ are the complex conjugates of $K_{11},K_{1\bar{1}}$,
respectively. We compute%
\begin{eqnarray}
K_{11,0} &=&T[{\theta ^{\bar{1}}(}KZ_{1})]-2{\omega _{1}}^{1}(T)K_{11}
\label{3.4} \\
\ &=&(L_{T}{\theta ^{\bar{1}})(}KZ_{1})+{\theta ^{\bar{1}}[(}L_{T}K)Z_{1}]+{%
\theta ^{\bar{1}}[}K(L_{T}Z_{1})]-2{\omega _{1}}^{1}(T)K_{11}.  \notag
\end{eqnarray}

It is easy to compute the first term using the (complex conjugate of)
structure equation (\ref{2.2}) and the third term using the formula $%
[T,Z_{1}]$ $=$ $-A_{11}Z_{\bar{1}}$ $+$ ${\omega _{1}}^{1}(T)Z_{1}$ (\cite%
{L1}). For the second term, if we take $K$ to be the initial $CR$ structure $%
J_{(0)},$then%
\begin{equation}
L_{T}K=L_{T}J_{(0)}=2A_{J_{(0)},\theta }=0  \label{3.5}
\end{equation}

\noindent by (\ref{2.3}) and the assumption. So altogether we obtain%
\begin{equation}
K_{11,0}=A_{11}(K_{1\bar{1}}-K_{\bar{1}1})=2A_{11}K_{1\bar{1}}.  \label{3.6}
\end{equation}

Note that $K^{2}=-I$ implies that $K_{1\bar{1}}=\pm i(1+|K_{11}|^{2})^{\frac{%
1}{2}}$ and $K_{\bar{1}1}=-K_{1\bar{1}.}$ It follows that%
\begin{equation}
K_{1\bar{1},0}=-A_{11}K_{\bar{1}\bar{1}}+A_{\bar{1}\bar{1}}K_{11}.
\label{3.7}
\end{equation}

\noindent Here the point is that both $K_{11,0}$ and $K_{1\bar{1},0}$ are
linear in $A_{11}$ and $A_{\bar{1}\bar{1}}$ with coefficients being
\textquotedblright 0th-order\textquotedblright\ in a (co)frame. Using (3.6),
(3.7), we can express $(F_{J}K)_{,011}$ as follows:%
\begin{eqnarray}
(F_{J}K)_{,011} &=&iK_{1\bar{1}}K_{11,\bar{1}\bar{1}011}+iK_{\bar{1}\bar{1}%
}K_{11,\bar{1}1011}  \label{3.8} \\
&&-iK_{\bar{1}1}K_{\bar{1}\bar{1},11011}-iK_{11}K_{\bar{1}\bar{1},1\bar{1}%
011}+l.w.t..  \notag
\end{eqnarray}

\noindent Here and hereafter $l.w.t.$ will mean a lower weight term in $%
A_{11},A_{\bar{1}\bar{1}}$ up to weight 3 with coefficients in $K_{1\bar{1}%
}, $ $K_{\bar{1}\bar{1}},$ $K_{\bar{1}1},$ $K_{11}$ and their covariant
derivatives up to weight 5. Note that $A_{11},A_{\bar{1}\bar{1}}$ are of
weight 2 in $K_{1\bar{1}},K_{\bar{1}\bar{1}},$ $K_{\bar{1}1},$ $K_{11}.$ The
first four terms on the right-hand side of (\ref{3.8}) contain the hightest
weight terms of weight 4 in $A_{11},A_{\bar{1}\bar{1}}$ in view of the
commutation relations (\ref{2.10}) and (\ref{3.6}), (\ref{3.7}) as will be
shown below. Using (\ref{2.10}) repeatedly and (\ref{3.6}), we compute%
\begin{eqnarray}
K_{11,\bar{1}\bar{1}011} &=&K_{11,\bar{1}0\bar{1}11}+l.w.t.=K_{11,0\bar{1}%
\bar{1}11}-2K_{11}A_{\bar{1}\bar{1},1\bar{1}11}+l.w.t.  \label{3.9} \\
&=&2K_{1\bar{1}}A_{11,\bar{1}\bar{1}11}-2K_{11}A_{\bar{1}\bar{1},1\bar{1}%
11}+l.w.t..  \notag
\end{eqnarray}

Similarly we obtain%
\begin{eqnarray}
K_{11,\bar{1}1011}\ &=&2K_{1\bar{1}}A_{11,\bar{1}111}-2K_{11}A_{\bar{1}\bar{1%
},1111}+l.w.t.  \label{3.10} \\
\ K_{\bar{1}\bar{1},11011} &=&2K_{\bar{1}1}A_{\bar{1}\bar{1},1111}-2K_{\bar{1%
}\bar{1}}A_{11,\bar{1}111}+l.w.t.  \notag \\
K_{\bar{1}\bar{1},1\bar{1}011} &=&2K_{\bar{1}1}A_{\bar{1}\bar{1},1\bar{1}%
11}-2K_{\bar{1}\bar{1}}A_{11,\bar{1}\bar{1}11}+l.w.t..  \notag
\end{eqnarray}

Substituting (\ref{3.9}), (\ref{3.10}) in (\ref{3.8}), we get, in view of (%
\ref{3.3}),%
\begin{equation}
(D_{J}F_{J}K)_{11,0}=-2iA_{11,\bar{1}\bar{1}11}+2iA_{\bar{1}\bar{1}%
,1111}+l.w.t..  \label{3.11}
\end{equation}

Now substituting (\ref{3.2}) and (\ref{3.11}) in (\ref{3.1}) gives%
\begin{equation}
\dot{A}_{11}=-\frac{1}{3}A_{11,\bar{1}\bar{1}11}+A_{11,00}+\frac{2i}{3}A_{11,%
\bar{1}10}+l.w.t.  \label{3.12}
\end{equation}

\noindent (note that the \textquotedblright bad\textquotedblright terms
cancel). Define $L_{\alpha }A_{11}=-A_{11,1\bar{1}}-A_{11,\bar{1}1}+i\alpha
A_{11,0}$ for a complex number $\alpha $. Let $L_{\alpha }^{\ast }$ be the
formal adjoint of $L_{\alpha }.$ It is a direct computation (cf. p.1257 in 
\cite{CL2}) that

\begin{eqnarray}
L_{\alpha }^{\ast }L_{\alpha }A_{11} &=&2(A_{11,11\bar{1}\bar{1}}+A_{11,\bar{%
1}\bar{1}11})-i(3+\alpha +\bar{\alpha}-|\alpha |^{2})A_{11,1\bar{1}0}
\label{3.13} \\
&&+i(3-\alpha -\bar{\alpha}-|\alpha |^{2})A_{11,\bar{1}10}+l.w.t..  \notag
\end{eqnarray}

Using the commutation relations (\ref{2.10}), we can easily obtain%
\begin{eqnarray}
A_{11,11\bar{1}\bar{1}} &=&A_{11,\bar{1}\bar{1}11}+2iA_{11,\bar{1}%
10}+2iA_{11,1\bar{1}0}+l.w.t.  \label{3.14} \\
A_{11,00} &=&-iA_{11,1\bar{1}0}+iA_{11,\bar{1}10}+l.w.t..  \notag
\end{eqnarray}

In view of (\ref{3.14}) and (\ref{3.13}), we can rewrite (\ref{3.12}) as
follows:%
\begin{equation}
\dot{A}_{11}=-\frac{1}{12}L_{\alpha }^{\ast }L_{\alpha }A_{11}+l.w.t.
\label{3.15}
\end{equation}

\noindent for $\alpha =4+i\sqrt{3}.$ Since $\alpha $ is not an odd integer, $%
L_{\alpha }$ and hence $L_{\alpha }^{\ast }L_{\alpha }$ (note $L_{\alpha
}^{\ast }=L_{\bar{\alpha}}$) are subelliptic (e.g. \cite{CL1}). Taking the
complex conjugate of (\ref{3.15}) gives a similar equation for $A_{\bar{1}%
\bar{1}}$ only with $\alpha $ replaced by $-\bar{\alpha}.$ On the other
hand, we observe that $A_{11}=0,A_{\bar{1}\bar{1}}=0$ for all (valid) time
is a solution to (\ref{3.15}) and its conjugate equation (note that $l.w.t.$
vanishes if $A_{11}$ and $A_{\bar{1}\bar{1}}$ vanish as remarked
previously). Therefore by the uniqueness of the solution to a (or system of)
subparabolic equation(s), we conclude that $A_{11}$ stays vanishing under
the flow (\ref{1.8}).

\endproof%

\bigskip

\proof
\textbf{(of Corollary C)} By the existence of a short-time solution to (\ref%
{1.8}) (which is $C^{k}$ smooth for any given large $k,$ see \cite{CL1}) and
Lemma 3.1, we can find $\tau _{1}$ $>$ $0$ such that $A_{11(t)}=0$ for $%
0\leq t<\tau _{1}.$ It follows from (\ref{2.5}) that%
\begin{equation}
Q_{11(t)}=\frac{1}{6}W_{,11(t)}.  \label{4.9}
\end{equation}

\noindent Here $W_{,11(t)}=(Z_{1(t)})^{2}W_{(t)}-{\omega _{1}}%
_{(t)}^{1}(Z_{1(t)})Z_{1(t)}W_{(t)}$ and $W_{(t)}$ is the Tanaka-Webster
curvature with respect to $J_{(t)}$ (and fixed $\theta $). Therefore $u$ $=$ 
$-\frac{1}{6}W_{(t)}$ $+$ $\frac{i}{12}F_{J_{(t)}}J_{(0)}$ is a solution to (%
\ref{1.5}) by (\ref{4.9}) for $0\leq t<\tau _{2}\leq \tau _{1}$ with $\tau
_{2}$ so small that $W_{(t)}>0$ or $W_{(t)}<0$ (hence $\func{Re}u$ $=$ $-%
\frac{1}{6}W_{(t)}$ $\neq $ $0).$ Since (\ref{1.5}) is equivalent to (\ref%
{1.4}), we conclude the result by Theorem A (which still holds true in $%
C^{k} $ category for large $k)$.

\endproof%

\bigskip

\section{Proof of Theorem A}

Let $J_{(t)}$ be a solution to (\ref{1.1}) for $0\leq t<\tau $ with given
initial $J_{(0)}$ being fillable. We are going to construct an almost
complex structure $\check{J}$ \ on $M\times \lbrack 0,\tau ),$ integrable on 
$M\times (0,\tau ).$

There is a canonical choice of the (unitary) frame $Z_{1(t)}$ with respect
to $J_{(t)}$ (\cite{CL1}). Write $Z_{1(t)}=\frac{1}{2}(e_{1(t)}-ie_{2(t)})$
where $e_{1(t)},e_{2(t)}\in \xi $ and $J_{(t)}e_{1(t)}=e_{2(t)}.$ Let $%
\{\theta ,e_{(t)}^{1},e_{(t)}^{2}\}$ be a coframe dual to $%
\{T,e_{1(t)},e_{2(t)}\}$ on $M.$ We will identify $M\times \{t\}$ with $M$
(hence $T(M\times \{t\})$ with $TM$). Now we define an almost complex
structure $\check{J}$ \ at each point in $M\times \{t\}$ as follows:%
\begin{equation*}
\check{J}\text{ }|_{\xi }=J_{(t)},\text{ }\check{J}T=-a\frac{\partial }{%
\partial t}+bT+a(\alpha e_{1(t)}+\beta e_{2(t)}).
\end{equation*}

\noindent Here $a,b,\alpha ,\beta $ are some real (smooth) functions of
space variable and $t$, and $a\neq 0$ (so $\check{J}$ $\frac{\partial }{%
\partial t}$ is completely determined from the above formulas and $\check{J}%
^{2}=-identity.$ Strictly speaking, $\alpha ,\beta $ depend on the choice of
frame while $a,b$ are global). It is easy to see that the coframe dual to $%
\{e_{1(t)},$ $e_{2(t)},$ $\frac{\partial }{\partial t}-(b/a)T-\alpha
e_{1(t)}-\beta e_{2(t)},$ $(1/a)T\}$ is $\{e_{(t)}^{1}+\alpha dt,$ $%
e_{(t)}^{2}+\beta dt,$ $dt,$ $a\theta +bdt\}.$ So the following complex
1-forms:%
\begin{equation}
\Theta ^{1}=(e_{(t)}^{1}+\alpha dt)+i(e_{(t)}^{2}+\beta dt)=\theta
_{(t)}^{1}+\gamma ^{1}dt,  \label{4.1}
\end{equation}

\begin{equation}
\eta =(a\theta +bdt)-idt=a\theta +(b-i)dt  \label{4.2}
\end{equation}

\noindent are type (1,0) forms with respect to $\check{J}$. Here $\gamma
^{1}=\alpha +i\beta $ is really the $Z_{1(t)}$ coefficient of the vector
field $\alpha e_{1(t)}+\beta e_{2(t)}.$ Let $\Lambda ^{p,q}$ denote the
space of type (p,q) forms. The integrability of $\check{J}$ is equivalent to 
$d\Lambda ^{1,0}\subset \Lambda ^{2,0}+\Lambda ^{1,1}$ or $\Lambda
^{2,0}\wedge d\Lambda ^{1,0}=0.$ In terms of \ $\Theta ^{1},\eta ,$ the
integrability conditions read as follows:%
\begin{equation}
\eta \wedge \ \Theta ^{1}\wedge d\eta =0,  \label{4.3}
\end{equation}

\begin{equation}
\eta \wedge \ \Theta ^{1}\wedge d\Theta ^{1}=0.  \label{4.4}
\end{equation}

Substituting (\ref{4.1}), (\ref{4.2}) into (\ref{4.3}) and making use of $%
d\theta =d_{M}\theta =i\theta _{(t)}^{1}\wedge \theta _{(t)}^{\bar{1}}($here 
$d_{M}$ denotes the exterior differentiation on $M$ and $d=d_{M}+dt\frac{%
\partial }{\partial t}$ on $M\times (0,\tau )),$ we obtain%
\begin{equation*}
0=\eta \wedge \ \Theta ^{1}\wedge d\eta =[ab_{,\bar{1}}-(b-i)a_{,\bar{1}%
}+ia^{2}\gamma ^{1}]\theta \wedge \theta _{(t)}^{1}\wedge \theta _{(t)}^{%
\bar{1}}\wedge dt.
\end{equation*}

\noindent Here $b_{,\bar{1}}=Z_{\bar{1}(t)}b,$ $a_{,\bar{1}}=Z_{\bar{1}%
(t)}a. $ Therefore (\ref{4.3}) is equivalent to the relation between $a,b$
and $\gamma ^{1}$ as shown below:%
\begin{equation}
\gamma ^{1}=ia^{-1}b_{,\bar{1}}-ia^{-2}(b-i)a_{,\bar{1}}.  \label{4.5}
\end{equation}

Next note that $d\theta _{(t)}^{1}=d_{M}\theta _{(t)}^{1}+dt\wedge \dot{%
\theta}_{(t)}^{1}$ and%
\begin{equation*}
\dot{\theta}_{(t)}^{1}=-iE_{\bar{1}\bar{1}(t)}\theta _{(t)}^{\bar{1}}
\end{equation*}

\noindent ($E_{\bar{1}\bar{1}(t)}$ is the $\bar{1}\bar{1}-$component of $%
E_{J_{(t)}}$with respect to $J_{(t)}$). So substituting (\ref{4.1}), (\ref%
{4.2}) into (\ref{4.4}) and making use of (\ref{2.2}) for $\theta _{(t)}^{1}$%
, we obtain%
\begin{eqnarray*}
0 &=&\eta \wedge \Theta ^{1}\wedge d\Theta ^{1} \\
&=&(aiE_{\bar{1}\bar{1}(t)}+(b-i)A_{\bar{1}\bar{1}(t)}+a\gamma _{,\bar{1}%
}^{1})\theta \wedge \theta _{(t)}^{1}\wedge \theta _{(t)}^{\bar{1}}\wedge dt.
\end{eqnarray*}

\noindent Here $A_{\bar{1}\bar{1}(t)}$ is the $\bar{1}\bar{1}-$component of
the torsion tensor with respect to $J_{(t)}$ and $\gamma _{,\bar{1}}^{1}$ $%
:= $ $Z_{\bar{1}(t)}\gamma ^{1}$ $+$ ${\omega _{1}}_{(t)}^{1}(Z_{\bar{1}%
(t)})\gamma ^{1}$ where ${\omega _{1}}_{(t)}^{1}$ is the pseudohermitian
connection form with respect to $\theta _{(t)}^{1}$. Therefore (\ref{4.4})
is equivalent to the following relation between $a,b$ and $\gamma _{,\bar{1}%
}^{1}:$%
\begin{equation}
\gamma _{,\bar{1}}^{1}=-iE_{\bar{1}\bar{1}(t)}-a^{-1}(b-i)A_{\bar{1}\bar{1}%
(t)}.  \label{4.6}
\end{equation}

Substituting (\ref{4.5}) into (\ref{4.6}) and letting $f$ $=$ $a^{-1}$ $\neq 
$ $0,g$ $=$ $-ba^{-1},u$ $=$ $f$ $+ig$, we obtain an equation for a complex
valued function $u:$%
\begin{equation}
u_{,11}+iuA_{11(t)}=iE_{11(t)}.  \label{4.7}
\end{equation}

In view of (\ref{2.9}), we can express (\ref{4.7}) in an intrinsic form:%
\begin{equation}
J_{(t)}\circ D_{J_{(t)}}f+D_{J_{(t)}}g=E_{J_{(t)}}.  \label{4.8}
\end{equation}

Recall that $D_{J_{(t)}}f=\frac{1}{2}L_{X_{f}}J_{(t)}$ (\cite{CL1}) in which 
$X_{f}=-fT$ $+$ $i(Z_{1(t)}f)Z_{\bar{1}(t)}$ $-i(Z_{\bar{1}(t)}f)Z_{1(t)}$
is the infinitesimal contact diffeomorphism induced by $f.$ So the image of $%
D_{J_{(t)}}$ describes the tangent space of the orbit of the symmetry group
acting on $J_{(t)}$ by the pullback (in this case, the contact
diffeomorphisms are our symmetries). Now by assuming (\ref{1.4}) (which is
the same as (\ref{4.8})) in Theorem A, we obtain that $\check{J}$ is
integrable on $M\times (0,\tau ).$

On the other hand, $(M,J_{(0)})$ bounds a complex surface $N$ by our
assumption that $J_{(0)}$ is fillable. So we have another almost complex
structure $\hat{J}$ on $M\times (-\delta ,0]$ induced from $N$, integrable
on $M\times (-\delta ,0),$ and restricting to $J_{(0)}$ on $(M,\xi ).$ Up to
a diffeomorphism from $M\times (-\delta _{1},0]$ to $M\times (-\delta
_{2},0],$ identity on $M\times \{0\}$ for $\delta _{1,}\delta _{2}$ perhaps
smaller than $\delta ,$ we can assume that $\check{J}$ and $\hat{J}$
coincide at $M\times \{0\}$ where they may not coincide up to $C^{k}$ for $%
k\geq 1,$ however. We want to find a local diffeomorphism $\Phi $ from a
neighborhood $U$ of a point in $M$ times $(-\delta _{1},0]$ to a similar set
so that $\Phi $ is an identity on $U\times \{0\},$ and $\Phi ^{\ast }\hat{J}$
coincides with $\check{J}$ up to $C^{k}$ for some large integer $k$ at $%
U\times \{0\}$. Let $x^{i},0\leq i\leq 3$ denote the coordinates of $U\times
(-\delta _{1},0]$ with $x^{0}$ being the time variable for $(-\delta
_{1},0]. $ Let $y^{i},0\leq i\leq 3$ denote the corresponding coordinates of
the image of $\Phi $ with $y^{0}$ being the time variable. If we express $%
\hat{J}_{1}=\Phi ^{\ast }\hat{J}$ $=\Phi _{\ast }^{-1}(\hat{J}\circ \Phi
)\Phi _{\ast }$\ in coordinates, we usually write%
\begin{equation*}
(\hat{J}_{1})_{m}^{l}=\hat{J}_{i}^{j}\frac{\partial y^{i}}{\partial x^{m}}%
\frac{\partial x^{l}}{\partial y^{j}}
\end{equation*}

\noindent for $\hat{J}_{1}=$\ $(\hat{J}_{1})_{m}^{l}dx^{m}\otimes \frac{%
\partial }{\partial x^{l}}$ and $\hat{J}=$ $\hat{J}_{i}^{j}dy^{i}\otimes 
\frac{\partial }{\partial y^{j}}.$ Let $\eta =\Phi _{\ast }^{-1}.$ Then $%
\eta ^{-1}$ has the expression $(\frac{\partial y^{i}}{\partial x^{m}}),$
the Jacobian matrix of $\Phi ,$ in coordinates. We require $\eta
^{-1}=identity$ at each point with $x^{0}=0$ where $\hat{J}$ coincides with $%
\check{J}$. Differentiating $\hat{J}_{1}=\Phi ^{\ast }\hat{J}$ $=\Phi _{\ast
}^{-1}(\hat{J}\circ \Phi )\Phi _{\ast }=\eta (\hat{J}\circ \Phi )\eta ^{-1}$%
(considered as a matrix equation with respect to the above-mentioned bases)
in $x^{0}$ at $x^{0}=0$, we obtain%
\begin{equation}
\hat{J}_{1}^{\prime }-\hat{J}^{\prime }=\eta ^{\prime }\hat{J}-\hat{J}\eta
^{\prime }.  \label{4.10}
\end{equation}

\noindent Here the prime of $\hat{J}^{\prime }$ means the $y^{0}$-derivative
at $y^{0}=0$ while the prime of $\hat{J}_{1}^{\prime }$ and $\eta ^{\prime }$%
means the $x^{0}-$derivative at $x^{0}=0$. Finding $\Phi $ such that $\hat{J}%
_{1}=\Phi ^{\ast }\hat{J}$ \ coincides with $\check{J}$ up to $C^{1}$ at $%
U\times \{0\}$ is reduced to solving the above equation (\ref{4.10}) for $%
\eta ^{\prime }$ with $\hat{J}_{1}^{\prime }=\check{J}^{\prime }.$ Here the
prime of $\check{J}^{\prime }$ means the $t$-derivative at $t=0.$ And this
can be done by simple linear algebra as follows. First note that $C=\check{J}%
^{\prime }-\hat{J}^{\prime }$ satisfies $\hat{J}C+C\hat{J}=0$ since $\check{J%
}=\hat{J}$ at $U\times \{0\}$ and both $\check{J}^{\prime }$ and $\hat{J}%
^{\prime }$ satisfies the same relation as $C$ does. With respect to a
suitable basis, $\hat{J}$ has a canonical matrix representation:%
\begin{equation*}
\left( 
\begin{array}{cccc}
0 & -1 & 0 & 0 \\ 
1 & 0 & 0 & 0 \\ 
0 & 0 & 0 & -1 \\ 
0 & 0 & 1 & 0%
\end{array}%
\right)
\end{equation*}

\noindent Then $C$ has the matrix form $\left( 
\begin{array}{cc}
C_{11} & C_{12} \\ 
C_{21} & C_{22}%
\end{array}%
\right) $ where each $C_{ij}$ is a $2\times 2$ matrix $\left( 
\begin{array}{cc}
a_{ij} & b_{ij} \\ 
b_{ij} & -a_{ij}%
\end{array}%
\right) .$ Now the solution $\eta ^{\prime }$ to (\ref{4.10}) has the matrix
form $\left( 
\begin{array}{cc}
\eta _{11}^{\prime } & \eta _{12}^{\prime } \\ 
\eta _{21}^{\prime } & \eta _{22}^{\prime }%
\end{array}%
\right) $ where each $\eta _{ij}^{\prime }$ is a $2\times 2$ matrix $\left( 
\begin{array}{cc}
u_{ij} & v_{ij} \\ 
w_{ij} & s_{ij}%
\end{array}%
\right) $ satisfying the relations: $v_{ij}+w_{ij}=-a_{ij},$ $%
u_{ij}-s_{ij}=b_{ij}.$ Once $\eta ^{\prime }$ is determined by equation (\ref%
{4.10}), it is easy to construct the \textquotedblright
local\textquotedblright diffeomorphism $\Phi _{1}$ such that the inverse
Jacobian and its $x^{0}$-derivative at $x^{0}=0$ of $\Phi _{1}$ is $\eta =$%
the identity and $\eta ^{\prime }$, resp. (we may need to shrink the time
interval $(-\delta _{1},0]$). So if we start with $\hat{J}_{1}=\Phi
_{1}^{\ast }\hat{J}$ instead of $\hat{J}$ and repeat the above procedure
looking for $\Phi _{2}$ so that $\hat{J}_{2}=\Phi _{2}^{\ast }\hat{J}_{1}$
coincides with $\check{J}$ at $U\times \{0\}$ up to $C^{2}$, we
differentiate $\hat{J}_{2}=\eta _{1}\hat{J}_{1}\eta _{1}^{-1}$ twice with
respect to $x^{0}$\bigskip\ at $x^{0}=0.$ Here $\eta _{1}$ denotes the
inverse Jacobian matrix of $\Phi _{2}$ (to be determined). Requiring $\hat{J}%
_{2}^{\prime }=\hat{J}_{1}^{\prime }$ and $\eta _{1}=identity$ (at $x^{0}=0$%
) implies $\eta _{1}^{\prime }=0.$ It then follows that $\eta _{1}^{\prime
\prime },$the second derivative of $\eta _{1}$ in $x^{0}$ at $x^{0}=0,$
satisfies a similar equation as in (\ref{4.10}):%
\begin{equation}
\eta _{1}^{\prime \prime }\hat{J}_{1}-\hat{J}_{1}\eta _{1}^{\prime \prime }=%
\hat{J}_{2}^{\prime \prime }-\hat{J}_{1}^{\prime \prime }.  \label{4.11}
\end{equation}

We can verify that the right-hand side anti-commutes with $\hat{J}_{1}$ as
follows: $(\hat{J}_{2}^{\prime \prime }-\hat{J}_{1}^{\prime \prime })\hat{J}%
_{1}+$\ $\hat{J}_{1}(\hat{J}_{2}^{\prime \prime }-\hat{J}_{1}^{\prime \prime
})=(\hat{J}_{2}^{\prime \prime }\hat{J}_{1}+\hat{J}_{1}\hat{J}_{2}^{\prime
\prime })-(\hat{J}_{1}^{\prime \prime }\hat{J}_{1}+\hat{J}_{1}\hat{J}%
_{1}^{\prime \prime })=-2(\hat{J}_{2}^{\prime })^{2}+2(\hat{J}_{1}^{\prime
})^{2}=0$ (here we have used $\hat{J}_{2}=\hat{J}_{1}$, $\hat{J}_{2}^{\prime
}=\hat{J}_{1}^{\prime }$ and $J^{\prime \prime }J+2(J^{\prime
})^{2}+JJ^{\prime \prime }=0$ for any almost complex structure $J$ by
differentiating $J^{2}=-I$ twice. So we can solve (\ref{4.11}) for $\eta
_{1}^{\prime \prime }$ with $\hat{J}_{2}^{\prime \prime }=\check{J}^{\prime
\prime }$ and hence find a $\Phi _{2}$ with the required properties as
before. In general, suppose we have found $\Phi _{n-1}$such that $\hat{J}%
_{n-1}=$ $\Phi _{n-1}^{\ast }\hat{J}_{n-2}=\eta _{n-2}\hat{J}_{n-2}\eta
_{n-2}^{-1}$ coincides with $\check{J}$ up to $C^{n-1}$ at $x^{0}=0.$ Then
by the similar procedure we can find $\Phi _{n}$ such that $\hat{J}_{n}=\Phi
_{n}^{\ast }\hat{J}_{n-1}=\eta _{n-1}\hat{J}_{n-1}\eta _{n-1}^{-1}$
coincides with $\check{J}$ up to $C^{n}$ at $x^{0}=0,$ and the $x^{0}$%
-derivatives of $\eta _{n-1}$ vanish up to the order $n-1.$ Furthermore the $%
n$-th $x^{0}$-derivative $\eta _{n-1}^{(n)}$ satisfies a similar equation as
in (\ref{4.10}) or (\ref{4.11}):%
\begin{equation}
\eta _{n-1}^{(n)}\hat{J}_{n-1}-\hat{J}_{n-1}\eta _{n-1}^{(n)}=\check{J}%
^{(n)}-\hat{J}_{n-1}^{(n)}.  \label{4.12}
\end{equation}

\noindent Here $\check{J}^{(n)}$ denotes the $n$-th $t$-derivative of $%
\check{J}$ at $t=0$ while $\hat{J}_{n-1}^{(n)}$ means the $n$-th $x^{0}$%
-derivative of $\hat{J}_{n-1}$ at $x^{0}=0$.

Now $\hat{J}_{n}$ defined on $U\times (-\delta _{n},0]$ and $\check{J}$ \
defined on $U\times \lbrack 0,\delta _{n})$ for a small $\delta _{n}>0$
together form a $C^{n}$ integrable almost complex structure on $U\times (-$ $%
\delta _{n},\delta _{n})$. Therefore $U\times (-$ $\delta _{n},\delta _{n})$
is a complex manifold for $n\geq 4$ by a theorem of Newlander-Nirenberg (%
\cite{NN}). Since $M$ is compact, we can cover it by a finite number of $%
U^{\prime }s$ and have corresponding $\delta _{n}^{\prime }s$. For each
point in the overlap of two $U^{\prime }s$ considered in $U\times \{0\}$, we
can find local coordinate maps from an open neighborhood $V$ contained in
the intersection of two associated $U\times (-$ $\delta _{n},\delta
_{n})^{\prime }s$ into $C^{2}$ so that the transition map $\psi $ on the
\textquotedblright concave\textquotedblright\ part corresponding to positive
\textquotedblright time variable\textquotedblright\ is holomorphic (note
that our $(M,J_{(0)})$ is a strongly pseudoconvex boundary of $N.$ Moreover,
on $V\cap \{M\times \lbrack 0,\tau )\},$ we have the \textquotedblright
same\textquotedblright\ integrable almost complex structure $\check{J}$
while on the intersection of $V$ and $U\times (-$ $\delta _{n},0)^{\prime
}s, $ we may have \textquotedblright different\textquotedblright\ $(\hat{J}%
_{n})^{\prime }s$). We then extend $\psi $ to the pseudoconvex part
holomorphically, and denote the extension map by $\tilde{\psi}$. Now glue $%
V\cap \{U\times (-$ $\delta _{n},0)\}$ (complex structure $\hat{J}_{n}$)
with $V\cap \{$another copy of $U\times (-$ $\delta _{n},0)\}$ (perhaps
different $\hat{J}_{n}$) through $\tilde{\psi}.$ In this way we can manage
to extend the complex structure $\check{J}$ across $M\times \{0\}$ to $%
M\times (-\delta ,0)$ (globally) for some small $\delta >0$. Finally the
identity (a $CR$ diffeomorphism) on $(M,J_{(0)})$ extends to a
biholomorphism $\rho $ between $M\times (-\delta ,0)$ (perhaps smaller $%
\delta $) and an open set in $N$ near $M$ (recall that $N$ is a complex
surface that $M$ bounds). Glue $M\times (-\delta ,t)$ $(t<\tau )$ and $N$ \
via $\rho $ to form a complex surface $N_{t}$ that ($M\times \{t\},J_{(t)})$
bounds. We have shown that $J_{(t)}$ is fillable for $0<t<\tau .$

\section{\textbf{Describing embeddable $CR$ structures}}

Let $M\subset \mathbb{C}^{2}$ be a closed strongly pseudoconvex real
hypersurface with the inclusion map $i_{M}$ $:$ $M\rightarrow \mathbb{C}^{2}$%
. Let $\mathfrak{E}_{m}$ denote the (tame Fr\'{e}chet) manifold of $%
C^{\infty }$ smooth embeddings from $M$ into $\mathbb{C}^{2},$ which are
isotopic to $i_{M}$ (see, e.g., \cite{H}). Define an equivalent relation
\textquotedblleft $\sim $", in $\mathfrak{E}_{m}$ by 
\begin{eqnarray*}
\varphi ,\psi \in \mathfrak{E}_{m},\varphi \sim \psi &\Longleftrightarrow
&\exists \ \hbox{a $CR$ diffeomorphism}\ \rho :\varphi (M)\rightarrow \psi
(M) \\
&&\hbox{
such that }\ \rho \circ \varphi =\psi .
\end{eqnarray*}

Often we use the notation $\varphi ^{\ast }J_{\mathbb{C}_{2}}$ to mean $%
\varphi _{\ast }^{-1}\circ J_{\mathbb{C}_{2}}\circ \varphi _{\ast }$
restricted to $\varphi _{\ast }^{-1}\xi _{\varphi }$, where $J_{\mathbb{C}%
^{2}} $ is the complex structure of $\mathbb{C}^{2}$ and $\xi _{\varphi }$ $%
:=$ $T(\varphi (M))$ $\cap $ $J_{\mathbb{C}^{2}}(T(\varphi (M)))$. So $%
\varphi ^{\ast }J_{\mathbb{C}^{2}}$ is an induced $CR$ structure compatible
with the contact bundle $\varphi _{\ast }^{-1}\xi _{\varphi }$ on $M$. In
this notation, for $\varphi ,\psi \in \mathfrak{E}_{m}$ 
\begin{equation*}
\varphi \sim \psi \Longleftrightarrow \varphi ^{\ast }J_{\mathbb{C}%
^{2}}=\psi ^{\ast }J_{\mathbb{C}^{2}}.
\end{equation*}

Now let $\mathfrak{J}_{e}$ denote the set of all $CR$ structures on $M$,
induced from $\mathfrak{E}_{m}$ (or say, embeddable and realized by a nearby
embedding). That is to say, 
\begin{equation*}
\mathfrak{J}_{e}=\{\varphi ^{\ast }J_{\mathbb{C}^{2}}:\varphi \in \mathfrak{E%
}_{m}\}.
\end{equation*}%
(Or more precisely, $\mathfrak{J}_{e}=\{(M,\varphi _{\ast }^{-1}\xi
_{\varphi },\varphi ^{\ast }J_{\mathbb{C}^{2}}):\varphi \in \mathfrak{E}%
_{m}\}$ to specify the possibly different contact structures) Obviously $%
[\varphi ]\in \mathfrak{E}_{m}/_{\sim }\mapsto \varphi ^{\ast }J_{\mathbb{C}%
^{2}}\in \mathfrak{J}_{e}$ is a one-one correspondence. Next let $\xi $ $=$ $%
\xi _{i_{M}}$ denote the standard contact structure on $M$ induced from the
inclusion map $i_{M}.$ Let $\mathfrak{E}_{c}\subset \mathfrak{E}_{m}$ denote
the set of all contact embeddings $\varphi $: $M\rightarrow {\mathbb{C}^{2},}
$ i.e., $\varphi $ is an embedding such that $\varphi _{\ast }\xi $ $=$ $\xi
_{\varphi }$. Also let $Diff(M)$ and $Cont(M)$ denote the groups of
diffeomorphisms and contact (w.r.t. $\xi )$ diffeomorphisms, resp.. Then by
a theorem of Gray (\cite{G} or \cite{H}): two close enough contact
structures are isotopically equivalent, the map: $Cont(M)\backslash 
\mathfrak{E}_{c}\rightarrow Diff(M)\backslash \mathfrak{E}_{m}$ induced from
the identity is a one-one correspondence. Denote by $\mathfrak{J}_{c}$, the
set of all $CR$ structures on $M$, induced from $\mathfrak{E}_{c}$. That is
to say,%
\begin{equation*}
\mathfrak{J}_{c}=\{\varphi ^{\ast }J_{\mathbb{C}^{2}}:\varphi \in \mathfrak{E%
}_{c}\}.
\end{equation*}%
\noindent Note that elements in $\mathfrak{J}_{c}$, are all compatible with
the contact bundle $\xi $. Easy to see that the map:

\begin{equation}
\lbrack \varphi ]\in \mathfrak{E}_{c}/_{\sim }\rightarrow \varphi ^{\ast }J_{%
\mathbb{C}^{2}}\in \mathfrak{J}_{c}  \label{5.1}
\end{equation}%
is a one-one correspondence too. Let $\mathfrak{J}\supset \mathfrak{J}_{e}$
denote the set of all $CR$ structures on $M,$ whose associated contact
stuctures are isotropic to $\xi .$ Let $\mathfrak{J}_{\xi }\subset \mathfrak{%
J}$ denote those which are compatible with $\xi $. It follows that the maps: 
$Cont(M)\backslash \mathfrak{J}_{c}\rightarrow Diff(M)\backslash \mathfrak{J}%
_{e}$ and $Cont(M)\backslash \mathfrak{J}_{\xi }\rightarrow
Diff(M)\backslash \mathfrak{J}$ induced by the identity are also one-one
correspondences. Thus in summary we have the following commutative diagram:%
\begin{equation}
\begin{array}{ccccc}
Diff(M)\backslash \mathfrak{E}_{m}/_{\sim } & \twoheadrightarrow & 
Diff(M)\backslash \mathfrak{J}_{e} & \subset & Diff(M)\backslash \mathfrak{J}
\\[6pt] 
\Uparrow &  & \Uparrow &  & \Uparrow \\[6pt] 
Cont(M)\backslash \mathfrak{E}_{c}/_{\sim } & \twoheadrightarrow & 
Cont(M)\backslash \mathfrak{J}_{c} & \subset & Cont(M)\backslash \mathfrak{J}%
_{\xi }%
\end{array}
\label{5.2}
\end{equation}%
where the maps indicated by $\twoheadrightarrow $ or $\Uparrow $ are all
one-one correspondences. We can put, say, $C^{\infty }$ topology on the
above spaces in (\ref{5.2}) so that all the maps are homeomorphisms.

\bigskip

\textbf{Remark 5.1.} After we mod out the symmetry group, we get the genuine
spaces $Diff(M)\backslash \mathfrak{J}$ $(Diff(M)\backslash \mathfrak{J}%
_{e}, $ resp.) of all $CR$ structures (all $CR$ structures on $M,$
embeddable in $\mathbb{C}^{2}$ and realized by an embedding isotopic to $%
i_{M},$ resp.) on $M$. According to (\ref{5.2}), we can reduce our spaces to
those in the contact category. So forgetting about symmetries, we should try
to parametrize the space $\mathfrak{E}_{c}/_{\sim }$ in such a way that $%
\mathfrak{J}_{c}$, becomes a submanifold of $\mathfrak{J}_{\xi }$. We will
look at this problem infinitesimally in the next section.

\textbf{Remark 5.2.} Every object in this section can be similarly defined
for general $n$ provided $M$ is a closed strongly pseudoconvex hypersurface
in $\mathbb{C}^{n+1}$. The diagram (\ref{5.2}) still holds for general
dimensions.

\bigskip

\section{\textbf{Parametrizing $\mathfrak{J}_{c}$ infinitesimally}}

In this section, we will observe what the tangent space $T_{J}\mathfrak{J}%
_{c}$ of $\mathfrak{J}_{c}$, at $J$ would be if we know that $\mathfrak{J}%
_{c}$ has a manifold structure. In view of (\ref{5.1}), we look at $%
\mathfrak{E}_{c}$ first. Suppose $\mathfrak{E}_{c}$ is a good space
(manifold, say). Can we describe the tangent space $T_{i_{M}}\mathfrak{E}%
_{c} $ of $\mathfrak{E}_{c}$, at the inclusion map $i_{M}$ quantitatively?
Since the results below hold for general $n$, we will discuss the general
case hereafter.

Let $M\supset \mathbb{C}^{n+1}$ be the boundary of a strongly pseudoconvex
bounded domain, defined by $\gamma =0$ with $d\gamma \not=0$ at $M$. Let $%
\theta =-i\partial \gamma $ (type (1,0) in $\mathbb{C}^{n+1}$; when
restricted on $M$, $\theta $ is a contact form). Choose type (1, 0)-form $%
\theta ^{\alpha }$, $\alpha =1,\ldots ,n$, near $M$ (locally) such that $%
\theta $, $\theta ^{\alpha }$, $\alpha =1,\ldots ,n$, are independent and%
\begin{equation}
d\theta =ih_{\alpha \bar{\beta}}\theta ^{\alpha }\wedge \theta ^{\bar{\beta}%
}+\rho \theta \wedge \bar{\theta}  \label{6.1}
\end{equation}%
\noindent (summation convention hereafter) at (points of) $M$ (see Section
7: Appendix for a proof). Let $i\zeta $, $Z_{\alpha }$, $\alpha =1,\ldots ,n$
be type (1,0) vector fields dual to $\theta $, $\theta ^{\alpha }$, $\alpha
=1,\ldots ,n$. With this $\theta $ on $M$ (note that $d\theta =ih_{\alpha 
\bar{\beta}}\theta ^{\alpha }\wedge \theta ^{\bar{\beta}}$ on $M$ since $%
\theta =\bar{\theta}$ on $TM$), we have the pseudohermitian geometry and can
talk about the torsion and covariant derivatives, etc. (\cite{L1}). Recall
that $C^{\infty }(M,\mathbb{C})$ denote the set of $C^{\infty }$ smooth
complex-valued functions on $M$. Let $T_{i_{M}}\mathfrak{E}_{c}$ denote the
space of all $\partial F_{t}/\partial t|_{t=0}$ where $F_{t}\in \mathfrak{E}%
_{c}:M\subset \mathbb{C}^{n+1}\rightarrow M_{t}\subset \mathbb{C}^{n+1}$ is
a family of contact diffeomorphisms (w.r.t. contact structures induced from $%
\mathbb{C}^{n+1})$ with $F_{t}$ $=$ $i_{M}$ at $t$ $=$ $0.$

\bigskip

\textbf{Lemma 6.1. }\textit{$T_{i_{M}}\mathfrak{E}_{c}$ $($complex version$)$
$=\{Y_{f}=i(f\zeta +f^{,\alpha }Z_{\alpha }):f\in C^{\infty }(M,\mathbb{C}%
)\} $, i.e., every \textquotedblleft external" contact vector field on $M$
with value in $T_{1,0}\mathbb{C}^{n+1}$ is of the form $Y_{f}$ and vice
versa.}

\bigskip

\proof
Take $Y$ $\in $ $T_{i_{M}}\mathfrak{E}_{c}.$ That is to say, $Y$ is an
\textquotedblleft external" contact vector field with value in $T_{1,0}%
\mathbb{C}^{n+1}$, i.e. (real version) $2\mathbb{R}eY=\partial
F_{t}/\partial t|_{t=0}$ where $F_{t}\in \mathfrak{E}_{c}:M\subset \mathbb{C}%
^{n+1}\rightarrow M_{t}\subset \mathbb{C}^{n+1}$ is a family of contact
diffeomorphisms with $F_{t}$ $=$ $i_{M}$ at $t$ $=$ $0.$ Let $\tilde{F}_{t}$
be a ($C^{\infty })$ smooth extension of $F_{t}$ to $\mathbb{C}^{n+1}$ so
that $\gamma _{t}=\gamma \circ \tilde{F}_{t}^{-1}$ is a defining function of 
$M_{t}$. Let $\theta _{t}=-i\partial \gamma _{t}$ be a contact form on $%
M_{t} $. Since $F_{t}$ is contact, we have%
\begin{equation}
\tilde{F}_{t}^{\ast }\theta _{t}=\lambda _{t}\theta  \label{6.2}
\end{equation}%
\noindent on $M$. Differentiating (\ref{6.2}) in $t$ at $t=0$ gives%
\begin{equation}
L_{2\mathbb{R}e\widetilde{Y}}\theta +\frac{\partial \theta _{t}}{\partial t}%
\Big|_{t=0}=\lambda \theta  \label{6.3}
\end{equation}%
\noindent where $\widetilde{Y}$ is an extension of $Y$ near $M$ and $\lambda
=\partial \lambda _{t}/\partial t|_{t=0}$. Now write $Y=if\zeta +g^{\alpha
}Z_{\alpha }$. It follows from the basic formula: $L_{x}=d\circ
i_{X}+i_{X}\circ d$ where $i_{X}$ is the interior product in the direction $%
X $ that%
\begin{equation}
L_{2\mathbb{R}e\widetilde{Y}}\theta =df+ih_{\alpha \bar{\beta}}\theta
^{\alpha }\wedge \theta ^{\bar{\beta}}(2\mathbb{R}eY,\cdot )  \label{6.4}
\end{equation}%
\noindent on $T_{1,0}(M)\oplus T_{0,1}(M)$ by (\ref{6.1}), where $T_{1,0}(M)$
($T_{0,1}(M),$ resp.) denotes the space of type (1,0) ((0,1), resp.) tangent
vectors. On the other hand, compute%
\begin{eqnarray}
\frac{\partial \theta _{t}}{\partial t}\Big|_{t=0} &=&-i\partial \Big(\frac{%
\partial \gamma _{t}}{\partial t}\Big|_{t=0}\Big)  \label{6.5} \\
&=&-\partial (d\gamma (-2\mathbb{R}eY))(\hbox{since}\ \frac{\partial \tilde{F%
}_{t}^{-1}}{\partial t}\Big|_{t=0}=-2\mathbb{R}eY)  \notag \\
&=&-\partial (f-\bar{f})\text{ \ }(\text{by }d\gamma =\partial \gamma +\bar{%
\partial}\gamma ).  \notag
\end{eqnarray}

Applying (\ref{6.3}) to $Z_{\bar{\alpha}}$ and using (\ref{6.4}) and (\ref%
{6.5}), we obtain 
\begin{equation*}
f_{,\bar{\alpha}}+ih_{\beta \bar{\alpha}}g^{\beta }=0.
\end{equation*}%
\noindent Raising the indices gives $g^{\alpha }=if^{,\alpha }$, as
expected. Conversely, given $Y_{f}$, we need to construct a one-parameter
family of contact embeddings $\Sigma _{t}\in \mathfrak{E}_{c}$, such that $%
\dot{\Sigma}_{t}|_{t=0}$ $:=$ $\frac{d}{dt}|_{t=0}\Sigma _{t}=2\mathbb{R}%
eY_{f}$. First we can find embeddings $\varphi _{t}\in \mathfrak{E}_{m}$ for 
$t$ small with%
\begin{equation}
\dot{\varphi}_{t}|_{t=0}=2\mathbb{R}eY_{iImf}\quad \hbox{and}\ \varphi
_{0}=i_{M}.  \label{6.6}
\end{equation}%
\noindent Second, by Gray's theorem, there exists a one-parameter family of
diffeomorphisms $\psi _{t}$ on $M$ with $\psi _{0}=Id$ such that $\varphi
_{t}\circ \psi _{t}\in \mathfrak{E}_{c}$ for all $t$. Since $\psi
_{t}=(\varphi _{t}\circ \psi _{t})-\varphi _{t}$ at $t=0$ and $(\varphi
_{t}\circ \psi _{t})^{\prime }=2\mathbb{R}eY_{g}$ for some $g$ from the
previous argument, hence%
\begin{equation}
\dot{\psi}_{t}|_{t=0}=2\mathbb{R}eY_{h}  \label{6.7}
\end{equation}%
\noindent for some $h$ by the linearity of $Y_{f}$ in $f$. Moreover, note
that $\dot{\psi}_{t}|_{t=0}$ is a tangent vector field on $M$ and $\func{Im}%
\zeta $ is tangent to $M$. Therefore $h$ must be a real function. By a
result in \cite{CL2} (for the smooth tame category), there exists a
one-parameter family of contact diffeomorphisms $\rho _{t}$ on $M$ with $%
\rho _{0}=Id$ such that%
\begin{equation}
\dot{\rho}_{t}|_{t=0}=2\mathbb{R}eY_{\mathbb{R}ef-h}.  \label{6.8}
\end{equation}%
Now set $\Sigma _{t}=\varphi _{t}\circ \psi _{t}\circ \rho _{t}\in \mathfrak{%
E}_{c}$. Easy to see that $\dot{\Sigma}_{t}|_{t=0}=2\mathbb{R}%
e(Y_{iImf}+Y_{h}+Y_{\mathbb{R}ef-h})=2\mathbb{R}eY_{f}$ by (\ref{6.6}), (\ref%
{6.7}), and (\ref{6.8}).

\endproof%

\bigskip

We remark that $Y_{f}$ can be characterized by conditions: $Y_{f}$ $%
\lrcorner $ $\theta $ $=$ $f,$ $Y_{f}$ $\lrcorner $ $d\theta $ $=$ $-\bar{%
\partial}_{b}f$ ($\func{mod}$ $\bar{\theta}).$ Define a second order
operator $\mathfrak{D}_{J}$ as follows:%
\begin{equation}
\mathfrak{D}_{J}f=(f_{,\bar{\alpha}}\text{ }^{\beta }-iA_{\bar{\alpha}}\text{
}^{\beta }f)\theta ^{\bar{\alpha}}\otimes Z_{\beta }  \label{6.9}
\end{equation}%
\noindent for $f\in C^{\infty }(M,\mathbb{C}),$ say, where $A_{\bar{\alpha}}$
$^{\beta }$ is the torsion (\cite{L1}). Let $T^{\prime }$ denote $T_{1,0}%
\mathbb{C}^{n+1}$ restricted on $M$. We often do not distinguish between
bundles and their sections. Define $\bar{\partial}_{b}:T^{\prime
}\rightarrow T_{0,1}^{\ast }(M)\otimes T^{\prime }$ by

\begin{equation*}
\bar{\partial}_{b}(f^{\ell }\frac{\partial }{\partial z^{\ell }})=(\bar{%
\partial}f^{\ell })\otimes \frac{\partial }{\partial z^{\ell }}
\end{equation*}%
\noindent for $f^{\ell }\in C^{\infty }(M,\mathbb{C})$, $\ell =1,\ldots ,n+1$%
. Easy to check that the above definition is independent of choices of local
coordinates. For $\bar{Z}\in T_{0,1}(M)$, $Y\in T^{\prime }$, we have the
formula%
\begin{equation}
(\bar{\partial}_{b}Y)(\bar{Z})=\pi _{1,0}([(\bar{Z}\tilde{)},\tilde{Y}])%
\text{ restricted to }M  \label{6.10}
\end{equation}

\noindent where $(\bar{Z}\tilde{)},\tilde{Y}$ denote local ($C^{\infty }$
smooth) extensions of $\bar{Z},$ $Y$ to $\mathbb{C}^{n+1}$ near $M,$ resp.
and $\pi _{1,0}$ denotes the projection to $T_{1,0}\mathbb{C}^{n+1}.$

\bigskip

\textbf{Lemma 6.2.} \textit{$\bar{\partial}_{b}Y_{f}=i\mathfrak{D}_{J}f$}.

\bigskip

\proof
Let $\tilde{f}$ be a local extension of $f$ to $\mathbb{C}^{n+1}$ near $M.$
From (\ref{6.10}) we compute%
\begin{eqnarray}
(\bar{\partial}_{b}Y_{f})(Z_{\bar{\alpha}}) &=&\pi _{1,0}([Z_{\bar{\alpha}%
},Y_{\tilde{f}}]\text{ restricted to }M  \label{6.12} \\
&=&\theta ([Z_{\bar{\alpha}},Y_{\tilde{f}}])i\xi +\theta ^{\beta }([Z_{\bar{%
\alpha}},Y_{\tilde{f}}])Z_{\beta }.  \notag
\end{eqnarray}%
\noindent On the other hand,%
\begin{eqnarray}
\theta ([Z_{\bar{\alpha}},Y_{\tilde{f}}]) &=&Z_{\bar{\alpha}}(\theta (Y_{%
\tilde{f}}))-Y_{\tilde{f}}(\theta (Z_{\bar{\alpha}}))-d\theta (Z_{\bar{\alpha%
}},Y_{\tilde{f}})  \label{6.13} \\
&=&\tilde{f}_{,\bar{\alpha}}-0-\tilde{f}_{,\bar{\alpha}}=0  \notag
\end{eqnarray}%
\noindent by (\ref{6.1}) (or (\ref{A.1}) in Section 7: Appendix). We also
have%
\begin{eqnarray}
\theta ^{\beta }([Z_{\bar{\alpha}},Y_{\tilde{f}}]) &=&Z_{\bar{\alpha}%
}(\theta ^{\beta }(Y_{\tilde{f}}))-Y_{\tilde{f}}(\theta ^{\beta }(Z_{\bar{%
\alpha}}))-d\theta ^{\beta }(Z_{\bar{\alpha}},Y_{\tilde{f}})  \label{6.14} \\
&=&if_{,\bar{\alpha}}\text{ }^{\beta }-0+A_{\bar{\alpha}}\text{ }^{\beta }f 
\notag
\end{eqnarray}%
\noindent on $M$ by (\ref{A.2}) in Section 7: Appendix. Substituting (\ref%
{6.13}) and (\ref{6.14}) into (\ref{6.12}), we obtain the desired formula.

\endproof%

\bigskip

Define a map $\chi _{f}:M\rightarrow \mathbb{C}^{n+1}$ by%
\begin{equation}
\chi _{f}(p)=p+Y_{f}(p)  \label{6.15}
\end{equation}%
\noindent for $p\in M$ considered as a position vector in $\mathbb{C}^{n+1}$%
. For $f$ small, $\chi _{f}$ is an embedding.

\bigskip

\textbf{Lemma 6.3.} \textit{The map : $f\in Ker\mathfrak{D}_{J}\rightarrow
\chi _{f}\in \{CR$ maps : ($M,\xi ,J)$ $\rightarrow $ $\mathbb{C}^{n+1}\}$
is a one-one correspondence. Note that for $f$ small, the set of $CR$ maps
can be replaced by the set of $CR$ embeddings near the inclusion map }$i_{M}$%
\textit{.}

\bigskip

\proof
It is clear that $\chi _{f}$ is $CR$ for $f\in Ker\mathfrak{D}_{J}$ by Lemma
6.2. Also obviously the map is injective since $Y_{f}=Y_{g}$ implies $f=g$.
For surjectivity, we observe that $Y(p)=\varphi (p)-p$ for a given $CR$ map $%
\varphi :M\rightarrow \mathbb{C}^{n+1}$ is an \textquotedblleft external" $%
CR $ (hence contact) vector field on $M$. By Lemma 6.1, $Y=Y_{f}$ for some $%
f\in C^{\infty }(M,\mathbb{C})$. Moreover, we have $i\mathfrak{D}_{J}f=\bar{%
\partial}_{b}Y_{f}=\bar{\partial}_{b}Y=0$. Hence $f\in Ker\mathfrak{D}_{j}$
and $\varphi =\chi _{f}$.

\endproof%

\bigskip

\proof%
\textbf{\ (of Theorem D) }Note that the equivalence class $\{\varphi \in 
\mathfrak{E}_{c}:\varphi \sim i_{M}\}$ of $i_{M}$ is equal to the set of $CR$
embeddings near (isotopic to) $i_{M}$. Hence the tangent space $T_{[i_{M}]}(%
\mathfrak{E}_{c}/_{\sim })$ of $\mathfrak{E}_{c}/_{\sim }$ at [$i_{M}$] is
expected to identify with $\{Y_{f}\}/\{Y_{f}:CR\}$ which is in one-one
correspondence to $C^{\infty }(M,\mathbb{C})/Ker\mathfrak{D}_{J}$ according
to Lemma 6.1 and Lemma 6.2.

The map $\varphi _{t}\in \mathfrak{E}_{c}\rightarrow \varphi _{t}^{\ast }J_{%
\mathbb{C}^{2}}\in \mathfrak{J}_{c}$ induces a tangential map: 
\begin{equation*}
2\mathbb{R}eY_{f}\ \hbox{(real version)}=\varphi _{t}|_{t=0}\in T_{I}%
\mathfrak{E}_{c}\rightarrow \frac{d}{dt}\Big|_{t=0}\varphi _{t}^{\ast }J_{%
\mathbb{C}^{2}}=L_{2\mathbb{R}eY_{f}}J_{\mathbb{C}^{2}}\in T_{J}\mathfrak{J}%
_{c}.
\end{equation*}%
\noindent A similar argument as in the proof of Lemma 3.5 in \cite{CL1}
shows that%
\begin{equation}
L_{Y_{f}}J_{\mathbb{C}^{2}}=2\mathfrak{D}_{J}f=-2i\bar{\partial}_{b}Y_{f}
\label{6.17}
\end{equation}%
\noindent (the last equality is due to Lemma 6.2.). Hence 
\begin{equation*}
L_{2\mathbb{R}eY_{f}}J_{\mathbb{C}^{2}}=4\mathbb{R}e(\mathfrak{D}_{J}f)=4Im(%
\bar{\partial}_{b}Y_{f})
\end{equation*}%
\noindent by (\ref{6.17}) (note that the operator $\mathfrak{D}_{J}$ and the
real operator $D_{J}$ (=$B_{J}^{\prime }$ in \cite{CL1}) is related as
below: 
\begin{equation*}
2\mathbb{R}e\mathfrak{D}_{J}f=D_{J}f_{r}+JD_{J}f_{c}
\end{equation*}%
\noindent for $f=f_{r}+if_{c})$. In summary, we have shown the commutative
diagram (\ref{1.9}).

\endproof%

\bigskip

\section{Appendix}

Let $M$ be a strongly pseudoconvex real hypersurface in $\mathbb{C}^{n+1}$.
Let $\gamma $ be a defining function of $M$, i.e., $\gamma =0$, $d\gamma
\not=0$ at $M$. Set $\theta =-i\partial \gamma $. Then

\bigskip

\textbf{Lemma A.1.}\textit{There exist $\theta ^{\alpha }$, $\alpha
=1,\ldots ,n$ of type $(1,0)$-forms in $\mathbb{C}^{n+1}$ locally near $M$
such that $\theta $ and $\theta ^{\alpha }$, $\alpha =1,\ldots ,n$ are
independent and 
\begin{equation}
d\theta =ih_{\alpha \bar{\beta}}\theta ^{\alpha }\wedge \theta ^{\bar{\beta}%
}+\rho \theta \wedge \bar{\theta}  \label{A.1}
\end{equation}%
at $($points of$)$ $M$ with some function $\rho $ and $(h_{\alpha \bar{\beta}%
})$ being positive hermitian. Furthermore, $h_{\alpha \bar{\beta}}$ can be
chosen to be $c_{\alpha }\delta _{\alpha \beta }$ with $c_{\alpha }$
positive constant on $M$.}

\bigskip

\proof
Since $d\theta =-i\bar{\partial}\partial \gamma $ is of type (1,1), we have 
\begin{equation*}
d\theta =ih_{\alpha \bar{\beta}}\theta ^{\alpha }\wedge \theta ^{\bar{\beta}%
}+g_{\alpha }\theta \wedge \theta ^{\bar{\alpha}}+h_{\alpha }\theta ^{\alpha
}\wedge \bar{\theta}+\rho \theta \wedge \bar{\theta}
\end{equation*}%
for any choice of $\theta ^{\alpha }$, $\alpha =1,\ldots ,n$ (such that $%
\theta $ and $\theta ^{\alpha }$'s are independent). On $TM$, $\theta =\bar{%
\theta}$. Hence ($h_{\alpha \bar{\beta}}$) is hermitian and positive due to
strongly pseudoconveacity of $M$ and $g_{\alpha }=-\bar{h}_{\alpha }$ on $M$%
. Choose $U_{\beta }^{\alpha }$, $v^{\alpha }$ in $\mathbb{C}^{n+1}$ near $M$
such that $\tilde{h}_{\alpha \bar{\beta}}U_{m}^{\alpha }U_{\bar{\ell}}^{\bar{%
\beta}}=h_{m\bar{\ell}}$ with $\tilde{h}_{\alpha \bar{\beta}}$ (positive
hermitian) prescribed, and $i\tilde{h}_{\alpha \bar{\beta}}v^{\alpha }U_{%
\bar{\ell}}^{\bar{\beta}}=g_{\ell }$. Let $\tilde{\theta}^{\alpha }=U_{\beta
}^{\alpha }\theta ^{\beta }+v^{\alpha }\theta $. Then it follows that 
\begin{align*}
d\theta & =i\tilde{h}_{\alpha \bar{\beta}}\tilde{\theta}^{\alpha }\wedge 
\tilde{\theta}^{\bar{\beta}}+(h_{\alpha }+\bar{g}_{\alpha })\theta ^{\alpha
}\wedge \bar{\theta}+\rho \theta \wedge \bar{\theta}\ \hbox{near}\ M \\
& =i\tilde{h}_{\alpha \bar{\beta}}\tilde{\theta}^{\alpha }\wedge \tilde{%
\theta}^{\bar{\beta}}+\rho \theta \wedge \bar{\theta}\ \hbox{at}\ M.
\end{align*}%
$\ $\hfill

\endproof%

\bigskip

From the above proof we also see that there exist $\theta ^{\alpha }$'s such
that 
\begin{equation}
d\theta =ih_{\alpha \bar{\beta}}\theta ^{\alpha }\wedge \theta ^{\bar{\beta}%
}+f_{\alpha }\theta ^{\alpha }\wedge \bar{\theta}+\rho \theta \wedge \bar{%
\theta}  \label{A.1.1}
\end{equation}
\noindent near $M$ with $h_{\alpha \bar{\beta}}=\delta _{\alpha \beta }$ and 
$f_{\alpha }=0$ on $M$.

\bigskip

\textbf{Lemma A.2. }\textit{With $\theta ^{\alpha }$, $\alpha =1,\ldots ,n$
satisfying (}$\ref{A.1.1})$\textit{, there exist connection forms $\omega
_{\beta }^{\alpha }$, torsion $A_{\bar{\beta}}^{\alpha }$ such that 
\begin{equation}
d\theta ^{\alpha }=\theta ^{\beta }\wedge \omega _{\beta }^{\alpha }+A_{\bar{%
\beta}}^{\alpha }\theta \wedge \theta ^{\bar{\beta}}+\lambda \theta \wedge 
\bar{\theta}  \label{A.2}
\end{equation}%
with $A_{\bar{\beta}}^{\alpha }=A_{\bar{\alpha}}^{\beta }$ near $M$ and $%
\omega _{\beta }^{\alpha }+\omega _{\bar{\alpha}}^{\bar{\beta}}=0$ at $M$.}

\bigskip

Note that last terms in (\ref{A.1}) and (\ref{A.2}) disappear if restricted
to $TM$ (where $\theta =\bar{\theta})$.

\bigskip

\proof
It is clear that near $M$, 
\begin{equation}
d\theta ^{\alpha }=\theta ^{\beta }\wedge \omega _{\beta }^{\alpha }+\theta
\wedge \tau ^{\alpha }+\lambda \theta \wedge \bar{\theta}  \label{A.3}
\end{equation}%
for certain 1-forms $\omega _{\beta }^{\alpha }$ and $\tau ^{\alpha }$ being
of type $(0,1)$ since $\{\theta ,\theta ^{\alpha }=1,\ldots ,n\}$ spans $%
T_{1,0}\mathbb{C}^{n+1}$. Let $\omega _{\gamma \bar{\beta}}=\omega _{\gamma
}^{\beta }$, $\omega _{\bar{\beta}\gamma }=\overline{(\omega _{\beta
}^{\gamma })}$, $\tau _{\alpha }=\overline{(\tau ^{\alpha })}$ and $\tau _{%
\bar{\alpha}}=\tau ^{\alpha }$ (note that $h_{\alpha \bar{\beta}}=\delta
_{\alpha \beta }$). Write $df_{\alpha }-f_{\beta }\omega _{\alpha }^{\beta
}=f_{\alpha ,\beta }\theta ^{\beta }+f_{\alpha ,\bar{\beta}}\theta ^{\bar{%
\beta}}(\mathop {\hbox{\rm mod}\,}$ $\theta ,\bar{\theta})$. Now we
differentiate (\ref{A.1.1}) using (\ref{A.3}) to get 
\begin{eqnarray}
0 &=&i(-\omega _{\alpha \bar{\beta}}-\omega _{\bar{\beta}\alpha }+E_{\alpha 
\bar{\beta}})\wedge \theta ^{\alpha }\wedge \theta ^{\bar{\beta}}+i\theta
\wedge (\tau _{\bar{\beta}}\wedge \theta ^{\bar{\beta}})  \label{A.4} \\
&&-i\bar{\theta}\wedge (\tau _{\alpha }\wedge \theta ^{\alpha }+if_{\alpha
\beta }\theta ^{\beta }\wedge \theta ^{\alpha })(\mathop {\hbox{\rm mod}\,}%
\text{ }\theta \wedge \bar{\theta})  \notag
\end{eqnarray}%
\noindent where the error term $E_{\alpha \bar{\beta}}=-\delta _{\alpha
\beta }f_{\ell }\theta ^{\ell }+(if_{\alpha }\bar{f}_{\beta }-\rho \delta
_{\alpha \beta })\theta +(if_{\alpha \bar{\beta}}+\rho \delta _{\alpha \beta
})\bar{\theta}$ satisfies the condition: 
\begin{equation}
E_{\alpha \bar{\beta}}=0\ \hbox{at}\ M.  \label{A.5}
\end{equation}

From (\ref{A.4}) and writing $\tau ^{\alpha }=A_{\bar{\beta}}^{\alpha
}\theta ^{\bar{\beta}}$, we can express 
\begin{equation}
-\omega _{\alpha \bar{\beta}}-\omega _{\bar{\beta}\alpha }+E_{\alpha \bar{%
\beta}}=A_{\alpha \bar{\beta}\ell }\theta ^{\ell }+B_{\alpha \bar{\beta}\bar{%
\ell}}\theta ^{\bar{\ell}}  \label{A.6}
\end{equation}%
\noindent with $A_{\alpha \bar{\beta}\ell }=A_{\ell \bar{\beta}\alpha }$, $%
B_{\alpha \bar{\beta}\bar{\ell}}=B_{\alpha \bar{\ell}\bar{\beta}}$ and
clearly $A_{\bar{\beta}}^{\alpha }=A_{\bar{\alpha}}^{\beta }$. We adjust $%
\omega _{\alpha \bar{\beta}}$ by $\tilde{\omega}_{\alpha \bar{\beta}}=\omega
_{\alpha \bar{\beta}}+A_{\alpha \bar{\beta}\ell }\theta ^{\ell }$. Easy to
verify that (\ref{A.3}) holds with $\omega _{\bar{\beta}}^{\alpha }$
replaced by $\tilde{\omega}_{\beta }^{\alpha }$. Compute 
\begin{equation}
\tilde{\omega}_{\alpha \bar{\beta}}+\tilde{\omega}_{\bar{\beta}\alpha }=(A_{%
\bar{\beta}\alpha \bar{\ell}}-B_{\alpha \bar{\beta}\bar{\ell}})\bar{\ell}%
+E_{\alpha \bar{\beta}}  \label{A.7}
\end{equation}%
\noindent by (\ref{A.6}), where $A_{\bar{\beta}\alpha \bar{\ell}}=\overline{%
(A_{\beta \bar{\alpha}\ell })}$. Note that from (\ref{A.6}), $A_{\bar{\beta}%
\alpha \bar{\ell}}=B_{\alpha \bar{\beta}\bar{\ell}}$ at $M$ by hermitian
symmetry of $\omega _{\alpha \bar{\beta}}$ and (\ref{A.5}). Therefore it
follows from (\ref{A.7}) that at $M$ 
\begin{equation*}
\tilde{\omega}_{\alpha \bar{\beta}}+\tilde{\omega}_{\bar{\beta}\alpha }=0.
\end{equation*}

\endproof%

\end{document}